\documentclass[11pt,a4paper]{article}
\usepackage[T1]{fontenc}
\usepackage{lmodern}
\usepackage[margin=27mm]{geometry}
\usepackage{amsmath,amssymb,amsthm,mathtools,booktabs,enumitem,needspace,tikz}
\usetikzlibrary{arrows.meta,positioning,calc}
\usepackage{algorithm,algpseudocode}
\usepackage[colorlinks=true,linkcolor=blue,citecolor=blue,urlcolor=blue]{hyperref}
\setlist{itemsep=2pt,topsep=4pt}
\newtheorem{theorem}{Theorem}[section]
\newtheorem{proposition}[theorem]{Proposition}
\newtheorem{lemma}[theorem]{Lemma}
\newtheorem{corollary}[theorem]{Corollary}
\theoremstyle{definition}
\newtheorem{definition}[theorem]{Definition}
\newtheorem{example}[theorem]{Example}
\theoremstyle{remark}
\newtheorem{remark}[theorem]{Remark}
\newcommand{\F}{\mathbb F_2}
\newcommand{\nul}{\operatorname{null}}
\newcommand{\rk}{\operatorname{rank}}
\newcommand{\diag}{\operatorname{diag}}
\newcommand{\im}{\operatorname{im}}
\newcommand{\QQ}{\mathcal Q}
\newcommand{\cC}{\mathcal C}
\newcommand{\cF}{\mathcal F}
\newcommand{\cG}{\mathcal G}
\newcommand{\Rroot}[2]{\mathcal R_{#2}(#1)}
\newcommand{\eps}{\varepsilon}
\tikzset{
  gv/.style={circle,draw,fill=white,minimum size=3.6pt,inner sep=0pt,line width=.45pt},
  ge/.style={draw=black!70,line width=.5pt},
  le/.style={draw=red!75!black,line width=.7pt},
  ne/.style={draw=blue!65!black,line width=.7pt},
  fgdot/.style={circle,fill=black!85,inner sep=0pt,minimum size=2.3pt},
  fgedge/.style={draw=black!75,line width=.4pt},
}
\title{One-Vertex Lifts of Looped Circle Graphs:\\Obstructions and Interlace Polynomials}
\author{Qingying Deng\thanks{Corresponding author. Email: qingying@xtu.edu.cn. ORCID: \href{https://orcid.org/0000-0003-3931-435X}{0000-0003-3931-435X}.}, \quad	Yexiang Yan \\
	\medskip
	School of Mathematics and Computational Science, Xiangtan University, P.R. China \\
}
\date{}
\begin{document}
\maketitle
\begin{abstract}
Let $M$ be the adjacency matrix of a looped simple graph over $\F$ and let $t=\diag M$.
We study the one-vertex lift with adjacency matrix
$\left(\begin{smallmatrix}M+tt^T&t\\t^T&0\end{smallmatrix}\right)$.
A looped circle graph is called stable if its lift is a circle graph.
We prove that stability is equivalent to the existence of a pseudo-orientable bouquet realizing the given looped circle graph.
We then obtain a finite excluded toggle-minor characterization by applying two rooted constructions to the known excluded pivot-minors for circle graphs.
An exact enumeration yields $22$ minimal toggle-equivalence classes: one of order $5$, twelve of order $6$, and nine of order $7$.
The enumeration is accompanied by independently checkable linear-algebra and minor certificates.
Finally, we prove $q(G^\dagger;y)=yq(G;y)$ for the vertex-nullity interlace polynomial of every looped simple graph, and derive quasi-tree counting and real determinant formulas.
\end{abstract}

\section{Introduction}
Circle graphs are the intersection graphs of chords of a circle: the vertices are chords and two vertices are adjacent when the corresponding chords cross.
Their representations as chord diagrams are not unique, and this flexibility is central to their theory.
Adding a binary sign to every chord gives a \emph{signed chord diagram}, whose looped intersection graph has a loop at a chord exactly when the sign is $1$.
The same data can be read as a \emph{ribbon bouquet}: one vertex disc, one edge band per chord, with a band twisted exactly when the sign is $1$ (Section~\ref{sec:bouquets}); for a bouquet $B$, the resulting looped intersection graph is denoted $I(B)$.
Conversely, choosing a point in each attachment interval of a bouquet produces a signed chord diagram.
Thus looped circle graphs are precisely the looped intersection graphs of ribbon bouquets, although a given looped circle graph may have many such realizations.

Two complementary lines of work frame the problem.
On the graph side, Bouchet characterized circle graphs by three excluded vertex-minors~\cite{Bouchet1994}, and Geelen and Oum sharpened this to a characterization by fifteen excluded pivot-minors~\cite{GO2009}.
On the ribbon side, Ding and Kim introduced pseudo-orientable bouquets and an adjustment construction that converts such a bouquet into an orientable one~\cite{DK2026}.
Pseudo-orientability is defined by a \emph{certificate}: the vertex boundary is split into two arcs so that every untwisted edge has both endpoints on the same arc and every twisted edge has one endpoint on each arc.
The adjustment reverses one arc, toggles the signs of the edges that meet both arcs, and adds one untwisted edge between the two division points.
The property is existential in nature: it belongs to a particular realization and is not determined by the intersection matrix alone (Example~3.15 of~\cite{DK2026}).

The starting point of this paper is a one-vertex lift that connects these two pictures.
Let $G$ be a looped simple graph with adjacency matrix $M$, and put $t=\diag M$.
The \emph{one-vertex lift} $G^{\dagger}$ is the simple graph with adjacency matrix
\[
L(M)=\begin{pmatrix}M+tt^T&t\\t^T&0\end{pmatrix}.
\]
The matrix is alternating, and the new vertex is adjacent precisely to the looped vertices of $G$.
A looped circle graph $G$ is \emph{stable} if $G^{\dagger}$ is a circle graph (Definition~\ref{def:lift}).
The construction has an inverse: for a simple graph $H$ rooted at $p$, the \emph{rooted inverse} $\Rroot Hp$ has adjacency matrix $K+tt^T$, where $\left(\begin{smallmatrix}K&t\\t^T&0\end{smallmatrix}\right)$ is the adjacency matrix of $H$ with the root last.
The two constructions are mutually inverse (Proposition~\ref{prop:inverse}).
More importantly, the lift commutes with principal pivot transforms: for the parity map $\alpha$ of~\eqref{eq:alpha},
\[
L(M*X)=L(M)*\alpha(X)
\]
whenever $X$ is feasible (Lemma~\ref{lem:commute}), and the companion identity $\det L(M)[\alpha(X)]=\det M[X]$ (Lemma~\ref{lem:parity}) shows that the lift preserves all principal-minor information.
These two identities are the technical engine of the paper: they reduce questions about ribbon realizations and toggle-minors to finite linear-algebraic questions about matrices of order $|V(G)|$ or $|V(G)|+1$.
The matrix lift is the one used by Ding and Kim~\cite{DK2026}. Here we interpret it as an operation on looped graphs and use its rooted inverse to project circle-graph obstructions back to the original graph.

\textbf{Stability and pseudo-orientability.}
Our first main result characterizes stability in geometric terms.
Theorem~\ref{thm:realization} states that a looped circle graph $G$ is stable if and only if there is a pseudo-orientable bouquet $B$ with $I(B)\cong G$.
Our study of this connection was motivated by Ding and Kim's work on pseudo-orientable ribbon graphs~\cite{DK2026}.
In one direction, a chord representation of $G^{\dagger}$ is turned into a certificate by deleting the chord of the new vertex and reversing one of the two arcs.
In the other direction, the certificate is used to perform the adjustment and to recover $G^{\dagger}$ as the intersection graph of an orientable bouquet.
The proof also identifies the matrix mechanism: the chords crossing the deleted chord are exactly the chords that become twisted, and reversing one arc toggles precisely the adjacencies between pairs of such chords.
The existential quantifier in the theorem is essential, as the same looped circle graph can have bouquets with different pseudo-orientability properties.
Figure~\ref{fig:certificate} follows both directions on the three-vertex example of Figure~\ref{fig:lift}.

\textbf{Excluded toggle-minors.}
Our second main result is a finite forbidden-minor description of the same class.
The relevant relation is \emph{toggle equivalence}: two looped graphs are toggle-equivalent when one is a principal pivot transform of the other, up to isomorphism, and a \emph{toggle-minor} is obtained by a toggle followed by an induced subgraph, up to isomorphism (Definition~\ref{def:toggle}).
Stable looped circle graphs are closed under toggle-minors (Proposition~\ref{prop:closure}), so it is natural to ask for the minimal non-stable classes.
Starting from the fifteen excluded pivot-minors $\cF$ of Geelen and Oum, we form two rooted families.
Type~(1) takes $\Rroot Fr$ for $F\in\cF$ and $r\in V(F)$.
Type~(2) takes $\Rroot{F_N^+}r$, where $F_N^+$ is obtained from $F$ by adding a new vertex $r$ adjacent to exactly $N$, for a nonempty set $N\subseteq V(F)$.
After retaining the looped circle graphs and merging toggle-equivalent graphs, the minimal classes form a family $\cG$ (Definition~\ref{def:candidates}).
Theorem~\ref{thm:forbidden} states that a looped circle graph is stable if and only if it has no toggle-minor isomorphic to a member of $\cG$, and that $\cG$ is a complete set of representatives of the minimal non-stable toggle-equivalence classes.
The proof has two ingredients: a completeness lemma (Lemma~\ref{lem:complete}) producing a candidate from an arbitrary circle-graph obstruction, using the rooted projection of Corollary~\ref{cor:projection}; and a finite descent (Lemma~\ref{lem:descent}) reducing every candidate to a minimal one.
The resulting enumeration is exact rather than heuristic.
Proposition~\ref{prop:enumeration} lists $22$ minimal classes: one of order $5$, twelve of order $6$, and nine of order $7$.
Each step is accompanied by checkable linear-algebra and minor certificates.
The programs, data, and certificates are provided as ancillary files accompanying this arXiv submission, in \texttt{anc/enumeration/}.

\textbf{Interlace polynomials and quasi-trees.}
The lift yields one further consequence, a polynomial identity with enumerative content.
For a looped simple graph $G$, the vertex-nullity interlace polynomial is
$q(G;y)=\sum_{X\subseteq V(G)}(y-1)^{\nul M(G)[X]}$ (Definition~\ref{def:q}); it is the vertex-nullity specialization of the interlace polynomial introduced in~\cite{ABS2004}.
The algebraic input is a nullity pair: for $K=M+tt^T$ and $L(M)$ as above, the unordered pair $\{\nul K,\nul L(M)\}$ equals $\{\nul M,\nul M+1\}$ (Lemma~\ref{lem:nullity}).
Pairing each subset $X\subseteq V(G)$ with $X\cup\{p\}$ then gives
$q(G^{\dagger};y)=y\,q(G;y)$ (Theorem~\ref{thm:polynomial}).
Two applications of this identity are worth singling out.
First, $q$ is invariant under toggle equivalence (Corollary~\ref{cor:toggle-invariance}).
Second, for a bouquet $B$ the spanning quasi-trees are counted by $|\QQ(B)|=q(I(B);1)$; for a stable graph this count admits a real determinant formula, namely $|\QQ(B)|=q(G;1)=\det_{\mathbb R}(I+R)$ for every bouquet $B$ with $I(B)\cong G$, where $R$ is a real principal-minor detector of $M$ (Corollary~\ref{cor:determinant}).
The pseudo-orientable adjustment preserves the number of spanning quasi-trees (Corollary~\ref{cor:quasitree}), which is the enumerative counterpart of Theorem~\ref{thm:realization}; the determinant formula is the lift counterpart of the orientable Matrix--Quasi-tree Theorem~\cite{MMN2025,DK2026}.

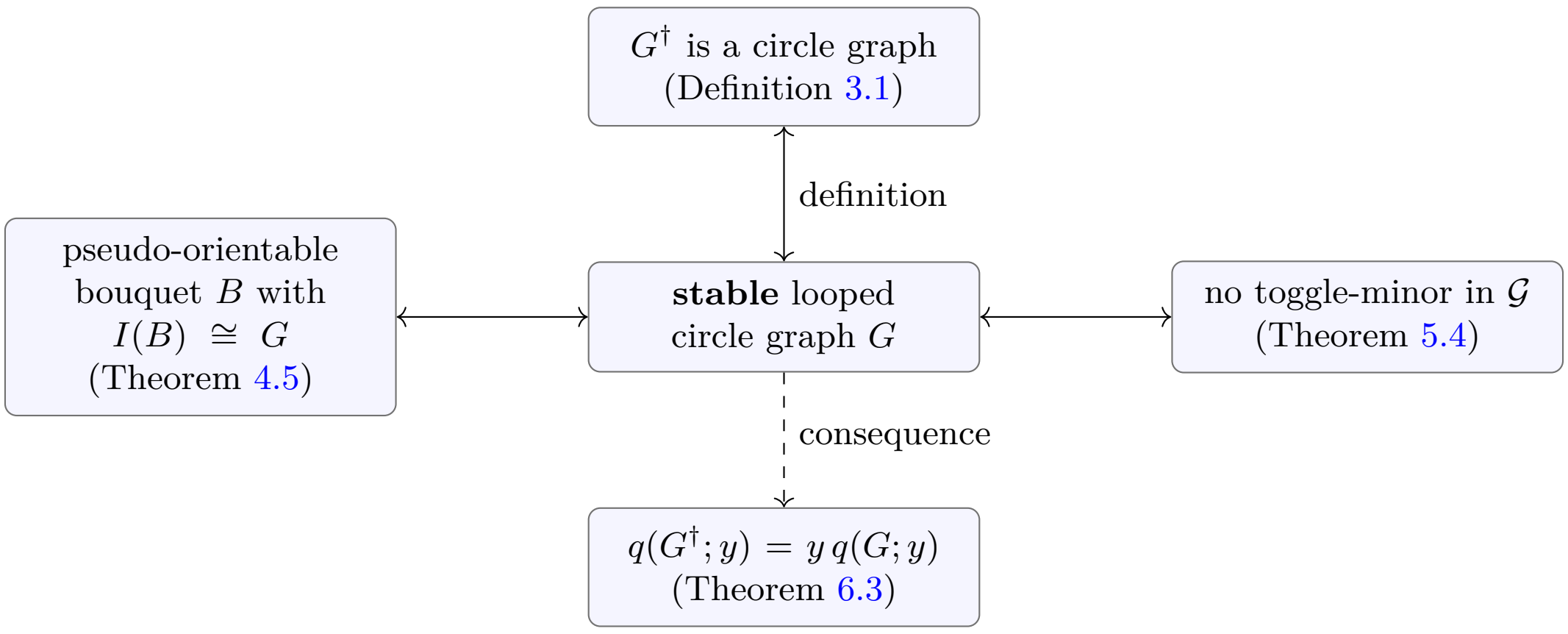
\begin{figure}[t]
\centering
\begin{tikzpicture}[every node/.style={font=\footnotesize},
  box/.style={draw=black!55,rounded corners=3pt,fill=blue!4,align=center,inner sep=4.5pt,text width=3.2cm}]
\node[box] (stable) at (0,0) {\textbf{stable} looped\\circle graph $G$};
\node[box] (lift) at (0,2.25) {$G^{\dagger}$ is a circle graph\\(Definition~\ref{def:lift})};
\node[box] (bouq) at (-5.25,0) {pseudo-orientable bouquet $B$ with $I(B)\cong G$\\(Theorem~\ref{thm:realization})};
\node[box] (minor) at (5.25,0) {no toggle-minor in $\mathcal G$\\(Theorem~\ref{thm:forbidden})};
\node[box] (poly) at (0,-2.25) {$q(G^{\dagger};y)=y\,q(G;y)$\\(Theorem~\ref{thm:polynomial})};
\draw[<->] (stable) -- (lift) node[midway,right] {definition};
\draw[<->] (stable) -- (bouq);
\draw[<->] (stable) -- (minor);
\draw[->,dashed] (stable) -- (poly) node[midway,right] {consequence};
\end{tikzpicture}
\caption{Three equivalent descriptions of stability. The upward arrow is the definition; the horizontal arrows are Theorem~\ref{thm:realization} and Theorem~\ref{thm:forbidden}. The interlace-polynomial identity holds for every looped simple graph, independently of stability.}
\label{fig:roadmap}
\end{figure}

The first two themes form the core of the paper; the third is a consequence of the same lift.
Section~\ref{sec:prelim} fixes notation, introduces chord diagrams and principal pivot transforms, and records the closure properties of looped circle graphs.
Section~\ref{sec:lift} develops the one-vertex lift, its rooted inverse, the parity correspondence, and the projection corollary with a fixed root.
Section~\ref{sec:bouquets} treats ribbon bouquets and proves Theorem~\ref{thm:realization}.
Section~\ref{sec:obstructions} introduces the two rooted constructions, proves the finite excluded-minor theorem, and describes the exact enumeration and its certificates.
Section~\ref{sec:polynomial} proves the nullity pair, the lift identity, and the determinant formula; it does not use the excluded-minor classification.
Appendix~\ref{app:representatives} displays one representative from each of the $22$ classes.
Figure~\ref{fig:roadmap} summarizes the logical relation among the three characterizations and the polynomial identity.
This paper is based on Chapters~3 and~4 of the second author's master's thesis~\cite{YanThesis}. The exposition and proofs have been revised for the present article.

\section{Graphs, matrices, and pivot operations}\label{sec:prelim}
\subsection{Conventions and intersection representations}
All sets, graphs, and matrices are finite.
Write $\F=\{0,1\}$ for the field with arithmetic modulo $2$.
For sets $X,Y$, their symmetric difference is $X\triangle Y=(X\setminus Y)\cup(Y\setminus X)$.
A square matrix indexed by $V$ has rows and columns both labeled by $V$; its restriction to rows and columns in $X$ is denoted by $M[X]$.
A \emph{principal minor} is $\det M[X]$, with $\det M[\varnothing]=1$.
We write $M^T$ for the transpose, $\diag M=(M_{vv})_{v\in V}$ for the diagonal column vector, and $I$ for an identity matrix of the required order.
The rank of a matrix is the dimension of its column space $\im M$, and its nullity $\nul M$ is the dimension of its kernel $\ker M$.
For an $n$-column matrix, $\nul M=n-\rk M$.
Unless a real matrix is explicitly specified, all matrix operations are over $\F$.
Whenever a binary scalar or vector occurs in a real identity, its integer $0/1$ representative is understood.
A binary matrix is \emph{symmetric} if $M^T=M$, and \emph{alternating} if it is symmetric and its diagonal is zero.
Equivalently, an alternating binary matrix satisfies $x^TMx=0$ for every column vector $x$.
A \emph{principal set} of a square matrix indexed by $V$ is a subset $X\subseteq V$; the corresponding \emph{principal submatrix} is $M[X]$.
Two square matrices indexed by the same set are \emph{matrix-isomorphic} if one is obtained from the other by a simultaneous permutation of rows and columns; we write $A\cong B$.

\begin{definition}[Looped graphs and induced subgraphs]\label{def:graphs}
A \emph{looped simple graph} is a pair $G=(V,E)$ in which
$E\subseteq\{\{v,w\}:v,w\in V\}$; a singleton $\{v\}\in E$ is a loop at $v$.
There are no multiple edges. A graph is \emph{simple} if it has no loops.
Its adjacency matrix $M(G)$ is given by $M(G)_{vw}=1$ if $\{v,w\}\in E$, and zero otherwise.
An \emph{isomorphism} is a bijection of vertex sets preserving adjacency, including loops; write $G\cong H$ when one exists.
The \emph{order} of $G$ is $|V|$.
For $S\subseteq V$, the \emph{induced subgraph} $G[S]$ has vertex set $S$ and all edges of $G$ contained in $S$.
Vertex deletion is $G-v=G[V\setminus\{v\}]$.
The open neighborhood is $N_G(v)=\{w\ne v:\{v,w\}\in E\}$.
A graph is connected if every two vertices can be joined by a sequence of edges between distinct consecutive vertices.
The \emph{underlying simple graph} $G^{-}$ is obtained by deleting all loops; thus $E(G^{-})=\{e\in E:|e|=2\}$.
\end{definition}

\begin{definition}[Chord diagrams and circle graphs]\label{def:chords}
A \emph{chord diagram} consists of a circle $C$, a set of $2n$ distinct points on $C$, and a partition of these points into $n$ labeled pairs, called chords.
Two chords \emph{cross} if their endpoints alternate in the cyclic order on $C$.
The \emph{intersection graph} has one vertex for each chord and one edge for each pair of crossing chords.
A \emph{circle graph} is a simple graph isomorphic to such an intersection graph.
A \emph{looped circle graph} is a looped simple graph whose underlying simple graph $G^{-}$ is a circle graph.

A \emph{signed chord diagram} additionally assigns a sign $\eps(e)\in\F$ to each chord $e$.
Its looped intersection graph has a loop at $e$ precisely when $\eps(e)=1$.
After a starting gap and a direction on $C$ are chosen, its cyclic order can be recorded by a \emph{double-occurrence word}: a word of length $2n$ in which every chord label appears twice.
Such a diagram or word \emph{realizes} a graph if its looped intersection graph is isomorphic to that graph.
\end{definition}

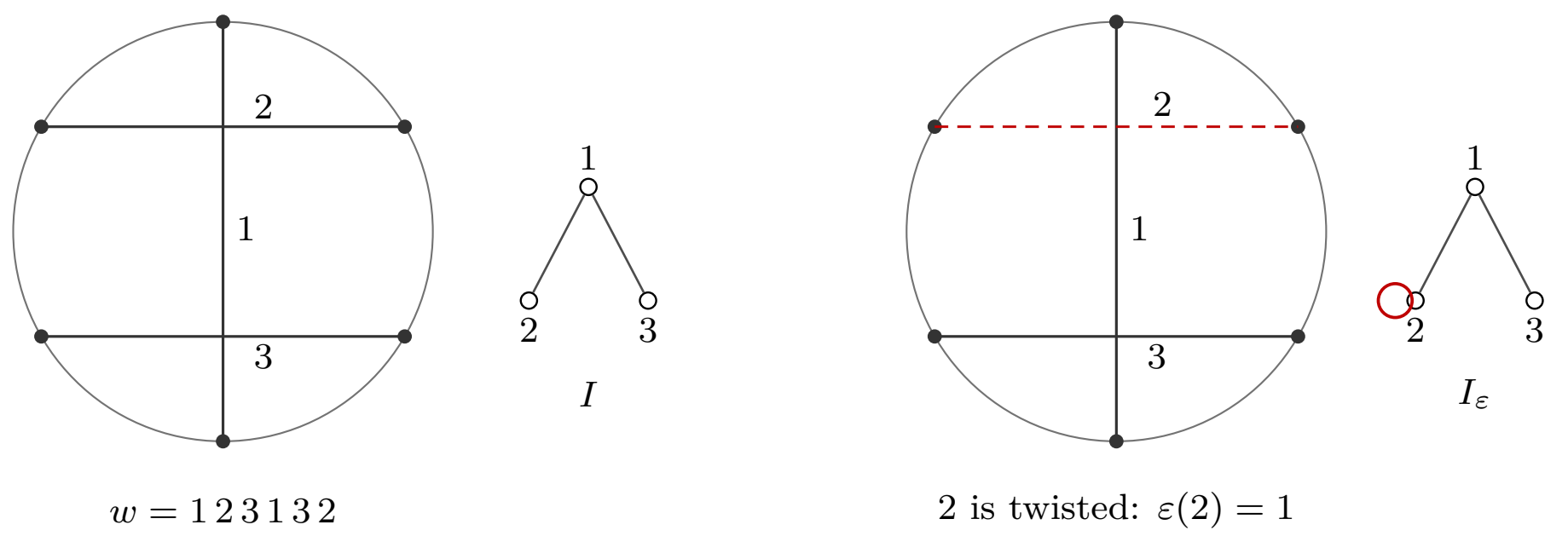
\begin{figure}[htbp]
	\centering
	\begin{tikzpicture}[scale=1.05,every node/.style={font=\scriptsize}]
		\begin{scope}
			\draw[black!55] (0,0) circle[radius=1.55];
			\foreach \a in {90,30,-30,-90,-150,-210}{
				\fill[black!80] (\a:1.55) circle[radius=1.5pt];
			}
			\draw[black!80,line width=.6pt] (90:1.55)--(-90:1.55);
			\draw[black!80,line width=.6pt] (30:1.55)--(-210:1.55);
			\draw[black!80,line width=.6pt] (-30:1.55)--(-150:1.55);
			\node[fill=white,inner sep=.8pt] at (0.17,0.02) {1};
			\node[fill=white,inner sep=.8pt] at (0.30,0.92) {2};
			\node[fill=white,inner sep=.8pt] at (0.30,-0.92) {3};
			\node at (0,-2.06) {$w=1\,2\,3\,1\,3\,2$};
		\end{scope}
		\begin{scope}[xshift=2.7cm,yshift=-0.15cm]
			\node[gv] (a1) at (0,0.48) {};
			\node[gv] (a2) at (-0.44,-0.36) {};
			\node[gv] (a3) at (0.44,-0.36) {};
			\draw[ge] (a1)--(a2); \draw[ge] (a1)--(a3);
			\node[above] at (a1) {1};
			\node[below] at (a2) {2};
			\node[below] at (a3) {3};
			\node at (0,-1.05) {$I$};
		\end{scope}
		\begin{scope}[xshift=6.6cm]
			\draw[black!55] (0,0) circle[radius=1.55];
			\foreach \a in {90,30,-30,-90,-150,-210}{
				\fill[black!80] (\a:1.55) circle[radius=1.5pt];
			}
			\draw[black!80,line width=.6pt] (90:1.55)--(-90:1.55);
			\draw[black!80,line width=.6pt] (-30:1.55)--(-150:1.55);
			\coordinate (A) at (-210:1.55);
			\coordinate (B) at (30:1.55);
			\draw[red!75!black,line width=.55pt,dash pattern=on 3pt off 2pt]
			(A)--(B);
			\node[fill=white,inner sep=.8pt] at (0.17,0.02) {1};
			\node[fill=white,inner sep=.8pt] at (0.34,0.94) {2};
			\node[fill=white,inner sep=.8pt] at (0.30,-0.92) {3};
			\node at (0,-2.06) {$2$ is twisted: $\eps(2)=1$};
		\end{scope}
		\begin{scope}[xshift=9.25cm,yshift=-0.15cm]
			\node[gv] (b1) at (0,0.48) {};
			\node[gv] (b2) at (-0.44,-0.36) {};
			\node[gv] (b3) at (0.44,-0.36) {};
			\draw[ge] (b1)--(b2); \draw[ge] (b1)--(b3);
			\draw[le] ($(b2)+(-.15,0)$) circle[radius=.125];
			\node[above] at (b1) {1};
			\node[below] at (b2) {2};
			\node[below] at (b3) {3};
			\node at (0,-1.05) {$I_\eps$};
		\end{scope}
	\end{tikzpicture}
	\caption{Chord diagrams and intersection graphs. (a) Chords $1$ and $2$, and chords $1$ and $3$, cross, while chords $2$ and $3$ do not; the endpoint order is the double-occurrence word $w=1\,2\,3\,1\,3\,2$. (b) The same diagram with chord $2$ twisted (shown dashed), recorded as a loop in the intersection graph $I_\eps$ at vertex $2$.}
	\label{fig:chords}
\end{figure}

Figure~\ref{fig:chords} illustrates a chord diagram, a signed chord diagram, and the two intersection graphs.

\subsection{Local complementation and principal pivots}
For a simple graph $G$ and a vertex $v$, the \emph{local complement} $G*v$ is obtained by reversing adjacency between every pair of distinct vertices in $N_G(v)$ and leaving all other adjacencies unchanged.
A \emph{vertex-minor} is an induced subgraph of a graph obtained by a sequence of local complements, up to isomorphism.
For an edge $uv$ in a simple graph, its \emph{edge pivot} is $G*uv=G*u*v*u$.
Two simple graphs are \emph{pivot-equivalent} if they are related by edge pivots and isomorphisms; a \emph{pivot-minor} is an induced subgraph of a pivot-equivalent graph.

\Needspace{14\baselineskip}
\begin{definition}[Feasibility and principal pivot transforms]\label{def:ppt}
Let $M$ be a symmetric binary matrix indexed by $V$.
A set $X\subseteq V$ is \emph{feasible} if $M[X]$ is invertible, or equivalently $\det M[X]=1$.
For such $X$, write, with $X$ first,
\[
M=\begin{pmatrix}P&Q\\R&S\end{pmatrix},\qquad P=M[X].
\]
The \emph{principal pivot transform} on $X$ is
\begin{equation}\label{eq:ppt}
M*X=\begin{pmatrix}P^{-1}&P^{-1}Q\\RP^{-1}&S+RP^{-1}Q\end{pmatrix}.
\end{equation}
For $X=\varnothing$, this means $M*\varnothing=M$.
\end{definition}

The transform remains symmetric and satisfies
\begin{equation}\label{eq:tucker}
\det((M*X)[Y])=\frac{\det M[X\triangle Y]}{\det M[X]}.
\end{equation}
It is an involution: $(M*X)*X=M$.
If the successive transforms are defined, then
$(M*X)*Y=M*(X\triangle Y)$~\cite{Tsatsomeros2000}.
For a simple graph and adjacent vertices $u,v$, the transform on $\{u,v\}$ is its edge pivot, with the same vertex labels.

At a looped vertex $v$, the singleton transform is called \emph{looped local complementation}.
It preserves $M_{vv}=1$ and the row and column incident with $v$, and replaces
\[
M_{ij}\ \longmapsto\ M_{ij}+M_{iv}M_{vj}\qquad(i,j\ne v).
\]
In particular, it toggles the loops at the neighbors of $v$ as well as the relevant off-diagonal adjacencies.
This operation differs from the simple-graph local complement defined above.
To state graph identities without ambiguity, write $\operatorname{lc}_v(G)$ for the ordinary local complement of a looped graph: it changes only off-diagonal adjacency among distinct neighbors and preserves every loop.
If $a_v=0$ and $a_i=M_{iv}$ for $i\ne v$, its matrix is
\[
M(\operatorname{lc}_v(G))=M+aa^T+\operatorname{Diag}(a),
\]
where $\operatorname{Diag}(a)$ is the diagonal matrix with diagonal vector $a$.
When $v$ is looped, write $\lambda_v(G)$ for the different operation with matrix $M*\{v\}$.

\begin{definition}[Toggle equivalence and toggle-minors]\label{def:toggle}
Looped graphs $G,H$ are \emph{toggle-equivalent} if, up to simultaneous relabeling of rows and columns,
$M(H)=M(G)*X$ for some feasible set $X$.
Their \emph{toggle-equivalence class} is the set of graphs related in this way.
A graph $H$ is a \emph{toggle-minor} of $G$, written $H\preceq G$, if
\[
M(H)\cong (M(G)*X)[S]
\]
for a feasible $X$ and a vertex subset $S$; matrix isomorphism means simultaneous permutation of rows and columns.
It is a \emph{proper toggle-minor} if $|V(H)|<|V(G)|$.
We use ``toggle'' specifically for this principal-pivot relation throughout the paper.
The \emph{toggle orbit} of $G$ is its toggle-equivalence class, that is, the set of all graphs $M(G)*X$ with $X$ feasible, taken up to isomorphism.
\end{definition}

\begin{lemma}[Factorization of feasible transforms]\label{lem:factor}
Let $M$ be a symmetric binary matrix and let $X$ be feasible for $M$.
Then $M*X$ is obtained from $M$ by a sequence of singleton transforms at looped vertices (indices $v$ with $M_{vv}=1$) and pair transforms at adjacent pairs of unlooped vertices (distinct $u,v$ with $M_{uu}=M_{vv}=0$ and $M_{uv}=1$).
\end{lemma}
\begin{proof}
Induct on $|X|$. For $X=\varnothing$ there is nothing to prove.
Suppose $X\ne\varnothing$.
If some $v\in X$ is looped, then $M_{vv}=1$ and $\{v\}$ is feasible.
By~\eqref{eq:tucker},
\[
\det\bigl((M*\{v\})[X\setminus\{v\}]\bigr)
=\frac{\det M[X]}{\det M[\{v\}]}=1,
\]
so $X\setminus\{v\}$ is feasible for $M*\{v\}$.
The composition law gives $M*X=(M*\{v\})*(X\setminus\{v\})$, and the induction hypothesis factors the second transform.
If no vertex of $X$ is looped, then $M[X]$ is alternating.
Since $\det M[X]=1$, the matrix $M[X]$ is nonzero, so there are distinct $u,v\in X$ with $M_{uv}=1$.
The principal block
$M[\{u,v\}]=\left(\begin{smallmatrix}0&1\\1&0\end{smallmatrix}\right)$
is invertible, so $\{u,v\}$ is feasible, and $u,v$ are adjacent unlooped vertices.
Applying~\eqref{eq:tucker} once more,
\[
\det\bigl((M*\{u,v\})[X\setminus\{u,v\}]\bigr)
=\frac{\det M[X]}{\det M[\{u,v\}]}=1,
\]
and the composition law together with induction finishes the proof.
\end{proof}

Looped circle graphs are closed under these transforms and under induced subgraphs~\cite{Chun2019,DK2026}.
In particular, all toggle-minors of a looped circle graph are looped circle graphs.

The relation $\preceq$ is transitive.
To see this, restriction to a set $S$ commutes with a transform on a subset $Y\subseteq S$ whenever that transform is defined.
Consequently, a sequence consisting of a transform on $X$, restriction to $S$, a transform on $Y\subseteq S$, and a further restriction to $T\subseteq S$ gives
$\bigl(M*(X\triangle Y)\bigr)[T]$.
The feasibility needed here follows from~\eqref{eq:tucker}.
Thus toggle-minors form a partial order on toggle-equivalence classes: two classes comparable in both directions must have the same order, no vertex deletion occurs, and they are the same class.

\section{The one-vertex lift and its rooted inverse}\label{sec:lift}
\begin{definition}[One-vertex lift and stability]\label{def:lift}
Let $G$ be a looped simple graph on $V$, let $M=M(G)$, and put $t=\diag M$.
For a new vertex $p\notin V$, define
\begin{equation}\label{eq:lift}
L(M)=\begin{pmatrix}M+tt^T&t\\t^T&0\end{pmatrix}.
\end{equation}
The simple graph with this adjacency matrix is the \emph{one-vertex lift} $G^{\dagger(p)}$, abbreviated $G^\dagger$.
A looped circle graph $G$ is \emph{stable} if $G^\dagger$ is a circle graph.
This does not depend on the name chosen for $p$.
\end{definition}

\begin{remark}[Relation to Ding--Kim]\label{rem:DK}
The matrix operation in~\eqref{eq:lift} is the binary specialization of the lift used by Ding and Kim~\cite[Lemmas 2.29 and 2.31]{DK2026}:
their matrix $\left(\begin{smallmatrix}A&v\\-v^T&0\end{smallmatrix}\right)$ equals $L(M)$ when $A=M+tt^T$ and $v=t$.
Here we study the simple graph with adjacency matrix $L(M)$.
Theorem~\ref{thm:realization} identifies the intersection graph of the adjustment of a pseudo-orientable bouquet $B$ with $I(B)^\dagger$.
Stability is a condition on the looped intersection graph, whereas pseudo-orientability is a property of a particular bouquet; the theorem relates stability to the existence of such a realization.
\end{remark}

\begin{figure}[htbp]
	\centering
	\begin{tikzpicture}[scale=1.1,every node/.style={font=\scriptsize}]
		\newcommand{\panelabase}{%
			\node[gv] (v1) at (0,1.0) {}; \node[gv] (v2) at (-.85,.1) {}; \node[gv] (v3) at (.85,.1) {};
			\draw[ge] (v1)--(v2);}
		\newcommand{\panelafull}{%
			\panelabase
			\draw[ge] (v1)--(v3);}
		\begin{scope}
			\panelafull
			\draw[le] ($(v1)+(90:.16)$) circle[radius=.12];
			\draw[le] ($(v3)+(0:.16)$) circle[radius=.12];
			\node[above left=-1pt] at (v1) {1}; \node[below left=-1pt] at (v2) {2}; \node[below right=-1pt] at (v3) {3};
			\node at (0,-1.55) {(a) $G$};
		\end{scope}
		\begin{scope}[xshift=3.45cm]
			\panelafull \node[gv] (p) at (0,-1.0) {};
			\draw[ne] (p)--(v1); \draw[ne] (p)--(v3);
			\draw[le] ($(v1)+(90:.16)$) circle[radius=.12];
			\draw[le] ($(v3)+(0:.16)$) circle[radius=.12];
			\draw[ne] ($(p)+(270:.16)$) circle[radius=.12];
			\node[above left=-1pt] at (v1) {1}; \node[below left=-1pt] at (v2) {2}; \node[below right=-1pt] at (v3) {3}; \node[left] at (p) {$p$};
			\node at (0,-1.55) {(b) add looped $p$};
		\end{scope}
		\begin{scope}[xshift=6.9cm]
			\panelabase \node[gv] (p) at (0,-1.0) {};
			\draw[ne] (p)--(v1); \draw[ne] (p)--(v3);
			\draw[ne] ($(p)+(270:.16)$) circle[radius=.12];
			\node[above left=-1pt] at (v1) {1}; \node[below left=-1pt] at (v2) {2}; \node[below right=-1pt] at (v3) {3}; \node[left] at (p) {$p$};
			\node at (0,-1.55) {(c) $\lambda_p$};
			\node[align=center] at (0,-2.15) {loops at $1,3$ and\\edge $13$ toggled off};
		\end{scope}
		\begin{scope}[xshift=10.35cm]
			\panelabase \node[gv] (p) at (0,-1.0) {};
			\draw[ne] (p)--(v1); \draw[ne] (p)--(v3);
			\node[above left=-1pt] at (v1) {1}; \node[below left=-1pt] at (v2) {2}; \node[below right=-1pt] at (v3) {3}; \node[left=0.22cm] at (p) {$p$ (root)};
			\node at (0,-1.55) {(d) $G^{\dagger}$};
		\end{scope}
		\draw[-{Latex[length=2mm]}] (1.25,-0.35)--(2.15,-0.35);
		\draw[-{Latex[length=2mm]}] (4.7,-0.35)--(5.6,-0.35);
		\draw[-{Latex[length=2mm]}] (8.15,-0.35)--(9.05,-0.35);
	\end{tikzpicture}
	\caption{The one-vertex lift of Definition~\ref{def:lift}. (a) A looped graph $G$ with loops at $1$ and $3$ and edges $12,13$. (b) Adjoin a looped vertex $p$ adjacent to every looped vertex. (c) Looped local complementation at $p$ keeps the edges incident with $p$ and its loop, and toggles off the loops at $1,3$ together with the edge $13$. (d) Removing the loop at $p$ gives $G^{\dagger}$, here a path on four vertices. The rooted inverse $\mathcal R(\cdot,p)$ reverses the whole procedure (Proposition~\ref{prop:inverse}). Red marks loops of $G$; blue marks the new vertex and its incident edges.}
	\label{fig:lift}
\end{figure}
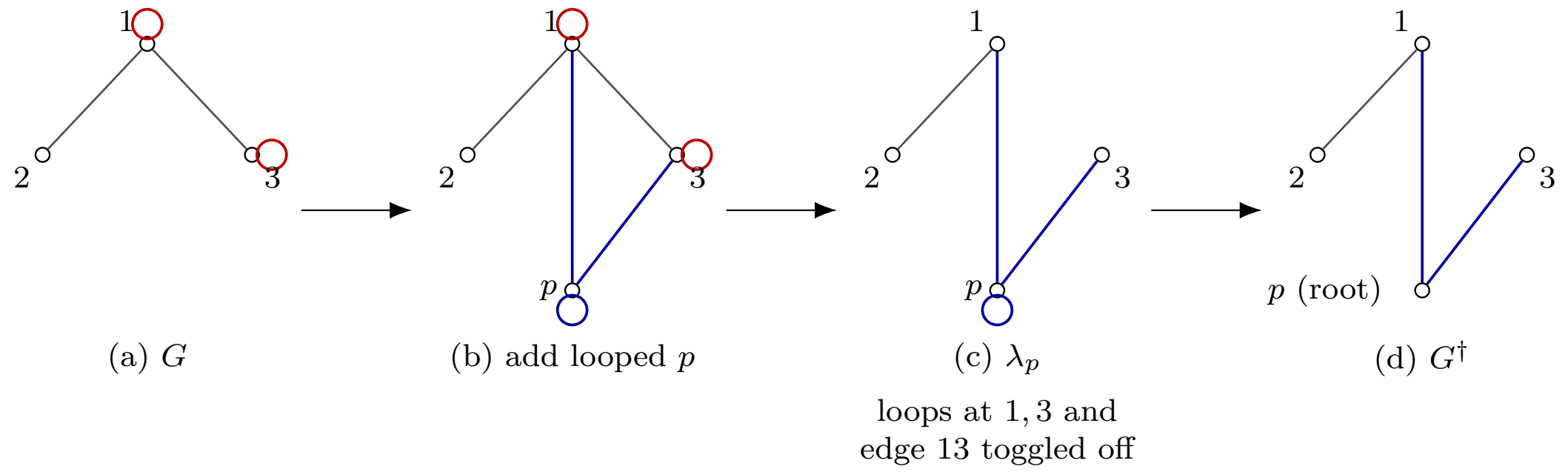

The matrix in~\eqref{eq:lift} is alternating, since $M_{ii}+t_i^2=0$ in $\F$.
Graphically, add a looped vertex adjacent to all looped vertices, apply looped local complementation there, and then remove the new vertex's loop.
Figure~\ref{fig:lift} shows the four steps on a three-vertex example.

\begin{definition}[Rooted inverse]\label{def:root}
A \emph{rooted simple graph} is a simple graph $H$ with a specified vertex $p$, its root.
A \emph{rooted isomorphism} is an isomorphism of simple graphs carrying the root of one graph to the root of the other.
Write the adjacency matrix of $(H,p)$, with the root last, as
\[
M(H)=\begin{pmatrix}K&t\\t^T&0\end{pmatrix}.
\]
The \emph{rooted inverse} $\Rroot Hp$ is the looped graph on $V(H)\setminus\{p\}$ with matrix
\begin{equation}\label{eq:root}
M(\Rroot Hp)=K+tt^T.
\end{equation}
\end{definition}

\begin{proposition}[Mutually inverse constructions]\label{prop:inverse}
For every looped graph $G$ and rooted simple graph $(H,p)$,
\[
\Rroot{G^{\dagger(p)}}p=G,\qquad (\Rroot Hp)^{\dagger(p)}=H.
\]
Thus lifting is a bijection between looped graphs and simple graphs with a specified new root, with labels retained or transported under rooted isomorphisms.
\end{proposition}
\begin{proof}
Write $M(H)=\left(\begin{smallmatrix}K&t\\t^T&0\end{smallmatrix}\right)$, with $K$ alternating.
Put $K'=K+tt^T$ and $t'=\diag K'$.
Since $K$ is alternating, its diagonal is zero, so $t'=t$.
Lifting $K'$ therefore gives
\[
L(K')=\begin{pmatrix}K'+t't'^T&t'\\t'^T&0\end{pmatrix}
=\begin{pmatrix}K+tt^T+tt^T&t\\t^T&0\end{pmatrix}
=\begin{pmatrix}K&t\\t^T&0\end{pmatrix}
=M(H),
\]
because $2tt^T=0$ in $\F$.
Thus $(\Rroot Hp)^{\dagger(p)}=H$.

Conversely, for $M=M(G)$ the lift has matrix $L(M)=\left(\begin{smallmatrix}K&t\\t^T&0\end{smallmatrix}\right)$ with $K=M+tt^T$ and $t=\diag M$.
Applying the rooted inverse gives
\[
K+tt^T=M+tt^T+tt^T=M,
\]
so $\Rroot{G^{\dagger(p)}}p=G$.
\end{proof}

For $X\subseteq V$, set
\begin{equation}\label{eq:alpha}
\alpha(X)=\begin{cases}X,&|X|\text{ even},\\X\cup\{p\},&|X|\text{ odd}.
\end{cases}
\end{equation}
Thus $\alpha$ is a bijection from the subsets of $V$ to the even-cardinality subsets of $V\cup\{p\}$, and preserves symmetric difference.

\begin{lemma}[Principal-minor correspondence]\label{lem:parity}
For every symmetric binary matrix $M$ and every $X$,
\begin{equation}\label{eq:parity}
\det L(M)[\alpha(X)]=\det M[X].
\end{equation}
Every odd-order principal minor of $L(M)$ is zero.
\end{lemma}
\begin{proof}
Put $K=M+tt^T$, which is alternating. Over $\F$ we have $M[X]=K[X]+t[X]t[X]^T$.
The bordered determinant identity, valid also for singular $K[X]$, gives
\begin{equation}\label{eq:bordered}
\det\bigl(K[X]+t[X]t[X]^T\bigr)=\det K[X]+
\det\begin{pmatrix}K[X]&t[X]\\t[X]^T&0\end{pmatrix};
\end{equation}
this is~\cite[Lemma 2.29]{DK2026}.

If $|X|$ is even, then $\alpha(X)=X$ and $L(M)[\alpha(X)]=L(M)[X]=K[X]$.
The bordered matrix in~\eqref{eq:bordered} is the principal submatrix $L(M)[X\cup\{p\}]$, so it is alternating, and it has odd order $|X|+1$; hence its determinant is zero.
Therefore~\eqref{eq:bordered} reduces to $\det M[X]=\det K[X]=\det L(M)[\alpha(X)]$.

If $|X|$ is odd, then $\alpha(X)=X\cup\{p\}$ and $L(M)[\alpha(X)]$ is exactly the bordered matrix in~\eqref{eq:bordered}.
Now $K[X]$ is alternating of odd order, so $\det K[X]=0$, and~\eqref{eq:bordered} again gives $\det M[X]=\det L(M)[\alpha(X)]$.

For $X=\varnothing$ both sides equal $1$ by the determinant convention.
Finally, $L(M)$ is alternating, and an alternating matrix has even rank; hence every odd-order principal minor of $L(M)$ vanishes.
\end{proof}

\begin{lemma}[Commutation with principal pivots]\label{lem:commute}
If $X$ is feasible for $M$, then
\begin{equation}\label{eq:commute}
L(M*X)=L(M)*\alpha(X),
\end{equation}
with the same labels and the same new vertex $p$.
\end{lemma}
\begin{proof}
The transform on the right is defined by Lemma~\ref{lem:parity}: $\alpha(X)$ is feasible for $L(M)$, because
$\det L(M)[\alpha(X)]=\det M[X]=1$.
For each $Y\subseteq V$, the principal minor of $L(M)*\alpha(X)$ indexed by $\alpha(Y)$ is, by~\eqref{eq:tucker} with pivot set $\alpha(X)$,
\[
\begin{aligned}
\det\bigl((L(M)*\alpha(X))[\alpha(Y)]\bigr)
&=\frac{\det L(M)[\alpha(X)\triangle\alpha(Y)]}{\det L(M)[\alpha(X)]}\\
&=\det L(M)[\alpha(X)\triangle\alpha(Y)]\\
&=\det M[X\triangle Y]
=\det(M*X)[Y].
\end{aligned}
\]
Here the second equality uses $\det L(M)[\alpha(X)]=1$, the third uses
$\alpha(X)\triangle\alpha(Y)=\alpha(X\triangle Y)$ and Lemma~\ref{lem:parity}, and the last is~\eqref{eq:tucker} for $M$.
Thus the two matrices have the same principal minor on every even-cardinality set $\alpha(Y)$.
Both matrices are alternating, since principal pivot transforms preserve alternation, so every odd-order principal minor of either matrix vanishes.
Finally, a symmetric binary matrix is determined by its principal minors of orders one and two, since
$a_{ij}=\det A[\{i,j\}]+a_{ii}a_{jj}$ for $i\ne j$.
The matrices are therefore equal.
\end{proof}

\begin{corollary}[Projection with a fixed root]\label{cor:projection}
If a simple graph $H$ is pivot-equivalent to $G^{\dagger(p)}$, transport the label $p$ along any relabelings.
Then $\Rroot Hp$ is toggle-equivalent to $G$.
\end{corollary}
\begin{proof}
After transporting labels, the composition law combines the pivot sequence into
$M(H)=L(M)*Y$ for one feasible set $Y$.
The matrix $L(M)$ is alternating, so every odd-order principal block is singular; hence its invertible principal block $L(M)[Y]$ has even order.
It follows that $Y=\alpha(X)$ for $X=Y\setminus\{p\}$, with $p\in Y$ exactly when $|X|$ is odd.
By Lemma~\ref{lem:parity},
\[
\det M[X]=\det L(M)[\alpha(X)]=\det L(M)[Y]=1,
\]
so $X$ is feasible for $M$.
Lemma~\ref{lem:commute} then gives $M(H)=L(M)*\alpha(X)=L(M*X)$.
Applying the rooted inverse at $p$ to both sides yields
$\Rroot Hp=M*X$, which is toggle-equivalent to $M=M(G)$.
Hence $\Rroot Hp$ is toggle-equivalent to $G$.
\end{proof}

\begin{proposition}[Compatibility with graph operations]\label{prop:operations}
With labels fixed and new vertex $p$, the following identities hold:
\begin{enumerate}[label=\textup{(\arabic*)}]
\item if $v$ is looped, $(\lambda_v(G))^{\dagger(p)}=G^{\dagger(p)}*pv$;
\item if $v$ is unlooped, $(\operatorname{lc}_v(G))^{\dagger(p)}=\operatorname{lc}_v(G^{\dagger(p)})$;
\item if $u,v$ are adjacent and unlooped, $(G*uv)^{\dagger(p)}=G^{\dagger(p)}*uv$;
\item for every $v\in V(G)$, $(G-v)^{\dagger(p)}=G^{\dagger(p)}-v$.
\end{enumerate}
Thus looped local complementation at $v$ corresponds to pivoting the root edge $pv$ in the lift. Ordinary local complementation at an unlooped vertex corresponds to ordinary local complementation at the same vertex in the lift. In the latter case $p$ is not adjacent to $v$, so its neighbourhood is unchanged. These are equalities of labelled graphs; no additional relabelling is needed.
\end{proposition}
\begin{proof}
For (1), the singleton $\{v\}$ is feasible and $\alpha(\{v\})=\{p,v\}$, so Lemma~\ref{lem:commute} applies; $pv$ is indeed an edge of the lift.
For (3), the feasible block on $\{u,v\}$ is $\left(\begin{smallmatrix}0&1\\1&0\end{smallmatrix}\right)$ and its parity image is the same set, so the same lemma applies.

For (2), use the adjacency vector $a$ defined above and set $\widetilde a=(a^T,0)^T$.
Since $t_v=0$, this is also the adjacency vector of $v$ in the lift.
Ordinary local complementation preserves $t$, and hence
\[
L(M(\operatorname{lc}_v(G)))
=L(M)+\widetilde a\widetilde a^T+\operatorname{Diag}(\widetilde a)
=M(\operatorname{lc}_v(G^{\dagger(p)})).
\]
Finally, (4) is the restriction identity $L(M[S])=L(M)[S\cup\{p\}]$ for $S=V(G)\setminus\{v\}$: deleting $v$ before lifting and deleting $v$ after lifting both amount to taking the principal submatrix on $S\cup\{p\}$.
\end{proof}

\begin{corollary}[Sequences of operations]\label{cor:operation-sequences}
If $G'$ is obtained from $G$ by looped local complementations and pivots on edges with unlooped endpoints, then their lifts are pivot-equivalent.
If ordinary local complements at unlooped vertices are also allowed, their lifts are related by ordinary local complements and edge pivots.
\end{corollary}
\begin{proof}
Apply Proposition~\ref{prop:operations} to each step; loop status and adjacency conditions are evaluated in the current graph at that step.
An induction on the number of operations then gives the two assertions.
\end{proof}

\begin{proposition}[Minor closure]\label{prop:closure}
Stable looped circle graphs are closed under toggle-minors and under ordinary local complementation at unlooped vertices.
\end{proposition}
\begin{proof}
Let $G$ be stable, put $M=M(G)$, and let $H$ be a toggle-minor of $G$.
Up to isomorphism, $M(H)=(M*X)[S]$ for a feasible set $X$ and a set $S\subseteq V(G)$.

By Lemma~\ref{lem:commute}, $L(M*X)=L(M)*\alpha(X)$.
The matrix $L(M)$ is alternating, so it has no looped vertices and every feasible principal block has even order.
Lemma~\ref{lem:factor} therefore factors the transform on $\alpha(X)$ into pair transforms at adjacent unlooped vertices, that is, into edge pivots of the simple graph $G^{\dagger}$.
Circle graphs are closed under edge pivots~\cite{Bouchet1994}, so $L(M*X)$ is a circle graph.
Thus $M*X$ is stable, and it is a looped circle graph by the closure statement in Section~\ref{sec:prelim}.

For vertex deletion, the restriction identity
\[
L\bigl((M*X)[S]\bigr)=L(M*X)[S\cup\{p\}]
\]
exhibits the lift of $(M*X)[S]$ as an induced subgraph of the circle graph $L(M*X)$.
Hence $(M*X)[S]$ is stable, and it is a looped circle graph because looped circle graphs are closed under induced subgraphs.
Since $H\cong(M*X)[S]$, $H$ is stable.

For ordinary local complementation at an unlooped vertex $v$, the underlying simple graph of $\operatorname{lc}_v(G)$ is $\operatorname{lc}_v(G^{-})$, which is a circle graph because circle graphs are closed under local complementation~\cite{Bouchet1994}.
By Proposition~\ref{prop:operations}(2), the lift of $\operatorname{lc}_v(G)$ is $\operatorname{lc}_v(G^{\dagger})$.
Since $G^{\dagger}$ is a circle graph and circle graphs are closed under local complementation, this lift is again a circle graph.
Thus $\operatorname{lc}_v(G)$ is stable.
\end{proof}

\section{Ribbon bouquets and pseudo-orientable realizations}\label{sec:bouquets}
\subsection{Ribbon graphs, signs, and certificates}
\begin{definition}[Ribbon graphs and quasi-trees]\label{def:ribbon}
A \emph{ribbon graph} is a compact surface with boundary presented as a union of finitely many vertex discs and rectangular edge bands.
The vertex discs have pairwise disjoint interiors, as do the edge bands.
Each edge band meets the vertex discs in exactly two disjoint boundary intervals, its attachment intervals; there are no other intersections, and all attachment intervals on each vertex disc are pairwise disjoint.
Corners may be smoothed without changing the surface's homeomorphism type.
The \emph{underlying graph} has one vertex per vertex disc and one edge per band, with incidences specified by the attachment intervals.
A ribbon graph is \emph{connected} if its underlying graph is connected.
Two ribbon graphs are identified when there is a homeomorphism carrying vertex discs to vertex discs and edge bands to edge bands, respecting any specified labels.

A \emph{boundary component} of a ribbon graph is a connected component of the boundary of its surface.
For $X\subseteq E(B)$, the \emph{spanning ribbon subgraph} $B|X$ retains every vertex disc and exactly the edge bands in $X$.
For a connected ribbon graph $B$, a \emph{spanning quasi-tree} is a set $X$ such that the surface $B|X$ has exactly one boundary component.
Write $\QQ(B)$ for the family of such edge sets.
A \emph{bouquet} is a ribbon graph with exactly one vertex disc.
An edge incident twice with the same vertex is a \emph{ribbon loop}; every edge of a bouquet is a ribbon loop.
\end{definition}

\Needspace{14\baselineskip}
\begin{definition}[Orientability and signed representations]\label{def:orientability}
A ribbon graph is \emph{orientable} if its surface admits a consistent orientation, that is, local choices of orientation agreeing on overlapping surface charts.
Otherwise it is \emph{non-orientable}.
For an edge $e$ of a bouquet, the surface consisting of the vertex disc and the band $e$ is either an annulus or a M\"obius band.
In the first case $e$ is \emph{untwisted}, and in the second it is \emph{twisted}.
Choose one point in each attachment interval; the resulting paired points on the vertex boundary form a chord diagram.
Give $e$ sign $\eps(e)=0$ when it is untwisted and sign $1$ when it is twisted.
Conversely, a signed chord diagram specifies a bouquet by attaching the indicated untwisted or twisted band at each pair of points.
An empty signed diagram represents a single vertex disc.
The \emph{looped intersection graph} $I(B)$ is the graph of this signed diagram.
\end{definition}

The choice of points inside attachment intervals does not affect cyclic order.
A bouquet is orientable exactly when all its edges are untwisted: an orientation of its vertex disc extends across precisely the untwisted bands.
An edge of a bouquet and a loop at a vertex of $I(B)$ are distinct objects; the latter records whether the former is twisted.

\begin{definition}[Pseudo-orientable bouquets and adjustment]\label{def:pseudo}
A bouquet $B$ with vertex boundary $C$ is \emph{pseudo-orientable} if there are two closed arcs $S_1,S_2$ with
$C=S_1\cup S_2$ and $S_1\cap S_2=\{P,Q\}$, where $P,Q$ are distinct points outside the attachment intervals, such that
\begin{enumerate}[label=\textup{(\roman*)}]
\item the two endpoints of every untwisted edge lie in the interior of the same arc;
\item every twisted edge has one endpoint in each arc's interior.
\end{enumerate}
The pair $(S_1,S_2)$ is a \emph{certificate}.

For a certificate, define the \emph{adjustment} $\widehat B$ by the following signed-diagram operation.
Reverse the order of endpoints on $S_2$ while keeping the order on $S_1$.
Toggle the sign of every edge having one endpoint in each arc.
Add a new edge $e$ with endpoints at $P,Q$ and sign zero.
Realize the resulting signed diagram as a bouquet.
\end{definition}

\begin{figure}[htbp]
	\centering
	\begin{tikzpicture}[scale=1.0,every node/.style={font=\scriptsize}]
		\newcommand{\ribbon}[3]{%
			\coordinate (RM) at ($(#1)!0.5!(#2)$);
			\coordinate (RP) at ($(#1)!0.5!(#2)!0.12cm!90:(#2)$);
			\coordinate (RA) at ($(#1)+(RP)-(RM)$); \coordinate (RB) at ($(#1)-(RP)+(RM)$);
			\coordinate (RC) at ($(#2)+(RP)-(RM)$); \coordinate (RD) at ($(#2)-(RP)+(RM)$);
			\draw[#3,line width=.55pt] (RA)--(RD);
			\draw[#3,line width=.55pt] (RB)--(RC);}
		\begin{scope}
			\draw[blue!45,line width=3pt] (90:1.6) arc[start angle=90,end angle=-45,radius=1.6];
			\draw[orange!75,line width=3pt] (-45:1.6) arc[start angle=-45,end angle=-270,radius=1.6];
			\draw[black!55] (0,0) circle[radius=1.6];
			\foreach \a in {90,45,0,-45,-90,-135,180,135}{\fill[black!80] (\a:1.6) circle[radius=1.5pt];}
			\node[above] at (90:1.85) {$P$};
			\node[below right] at (-45:1.85) {$Q$};
			\node[blue!55] at (60:1.95) {$S_1$};
			\node[orange!85!black] at (-150:1.95) {$S_2$};
			\draw[red!75!black,line width=.55pt,dash pattern=on 3pt off 2pt] (45:1.6)--(-135:1.6);
			\draw[black!80,line width=.6pt] (-90:1.6)--(180:1.6);
			\draw[red!75!black,line width=.55pt,dash pattern=on 3pt off 2pt] (0:1.6)--(135:1.6);
			\node[red!75!black,fill=white,inner sep=.8pt] at (-0.22,0.02) {1};
			\node[fill=white,inner sep=.8pt] at (-0.70,-0.92) {2};
			\node[red!75!black,fill=white,inner sep=.8pt] at (0.55,0.75) {3};
			\node at (0,-2.1) {(a) $B$ with certificate $(S_1,S_2)$};
		\end{scope}
		\begin{scope}[xshift=6.6cm]
			\draw[black!55] (0,0) circle[radius=1.6];
			\foreach \a in {90,45,0,-45,-90,-135,180,135}{\fill[black!80] (\a:1.6) circle[radius=1.5pt];}
			\node[above] at (90:1.85) {$P$};
			\node[below right] at (-45:1.85) {$Q$};
			\draw[black!80,line width=.6pt] (45:1.6)--(180:1.6);
			\draw[black!80,line width=.6pt] (-135:1.6)--(135:1.6);
			\draw[black!80,line width=.6pt] (0:1.6)--(-90:1.6);
			\draw[blue!65!black,line width=1.2pt] (90:1.6)--(-45:1.6);
			\node[fill=white,inner sep=.8pt] at (-0.30,0.66) {3};
			\node[fill=white,inner sep=.8pt] at (-1.30,0.02) {2};
			\node[fill=white,inner sep=.8pt] at (0.92,-0.78) {1};
			\node[blue!65!black,fill=white,inner sep=.8pt] at (0.60,0.28) {$e$};
			\node at (0,-2.1) {(b) adjusted bouquet $\widehat B$};
		\end{scope}
		\begin{scope}[xshift=3.35cm,yshift=-0.1cm]
			\node[gv] (h1) at (0,0.45) {}; \node[gv] (h2) at (-0.4,-0.35) {}; \node[gv] (h3) at (0.4,-0.35) {};
			\draw[ge] (h1)--(h2); \draw[ge] (h1)--(h3);
			\draw[le] ($(h1)+(90:.15)$) circle[radius=.12];
			\draw[le] ($(h3)+(0:.15)$) circle[radius=.12];
			\node[above] at (h1) {1}; \node[below] at (h2) {2}; \node[below] at (h3) {3};
			\node at (0,-1.05) {$G$};
		\end{scope}
		\begin{scope}[xshift=9.95cm,yshift=-0.1cm]
			\node[gv] (k1) at (0,0.45) {}; \node[gv] (k2) at (-0.4,-0.35) {}; \node[gv] (k3) at (0.4,-0.35) {}; \node[gv] (k4) at (0.85,0.35) {};
			\draw[ge] (k1)--(k2); \draw[ge] (k1)--(k4); \draw[ge] (k4)--(k3);
			\node[above] at (k1) {1}; \node[below] at (k2) {2}; \node[below] at (k3) {3}; \node[right] at (k4) {$p$};
			\node at (0.2,-1.05) {$G^{\dagger}$};
		\end{scope}
	\end{tikzpicture}
	\caption{A certificate and the adjustment. (a) A bouquet $B$ realizing $G$, with certificate arcs $S_1$ (blue) and $S_2$ (orange) meeting at $P,Q$. The untwisted chord $2$ has both endpoints in $S_2$, while each twisted chord $1,3$ (shown dashed) has one endpoint in each arc. (b) Reversing $S_2$, toggling the signs of the two spanning chords, and adding the edge $e=PQ$ gives the orientable bouquet $\widehat B$. Its intersection graph is $G^{\dagger}$; the intersection graph of $B$ is $G$. These are the graphs of Figure~\ref{fig:lift}.}
	\label{fig:certificate}
\end{figure}
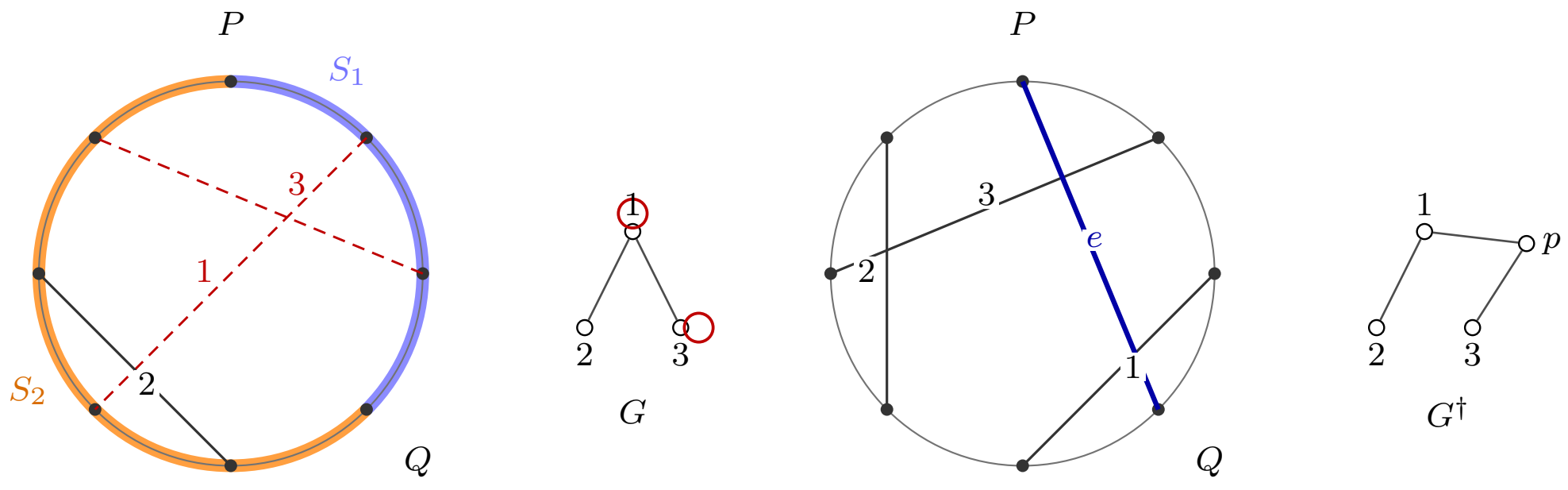

This combinatorial specification is the reversal and half-twist adjustment of~\cite{DK2026}.
Every edge of $\widehat B$ has sign zero, so it is orientable.
An orientable bouquet has a certificate: choose $P,Q$ in one gap with no endpoint between them on one of the two arcs.
Figure~\ref{fig:certificate} shows a certificate and the adjustment on the example of Figure~\ref{fig:lift}.

\begin{lemma}[Effect of reversing one arc]\label{lem:reversal}
Reversing the order on one of two complementary arcs changes the crossing relation of two chords if and only if each chord has one endpoint on each arc.
\end{lemma}
\begin{proof}
For two chords that both span the two arcs, alternation is determined by the relative orders of their two endpoints on the two arcs. Reversing exactly one order reverses alternation.
If both chords lie in one arc, reversing the order of all four endpoints preserves alternation.
If they lie in different arcs, they do not cross before or after reversal.
In the remaining case one chord spans the arcs and the other lies within one arc; they cross exactly when the endpoint of the spanning chord in that arc lies between the other chord's endpoints. Reversing that arc preserves betweenness.
\end{proof}

\begin{theorem}[Pseudo-orientable realizations]\label{thm:realization}
A looped circle graph $G$ is stable if and only if there exists a pseudo-orientable bouquet $B$ with $I(B)\cong G$.
\end{theorem}
\begin{proof}
Suppose first that $G$ is stable. Choose a chord representation of the circle graph $G^{\dagger(p)}$ and assign sign $0$ to every chord.
Let $e$ be the chord representing $p$, let $P,Q$ be its endpoints, and let $S_1,S_2$ be the two arcs into which $P,Q$ divide the boundary circle.
Write $M=M(G)$, $t=\diag M$, and $K=M+tt^T$.

In the matrix $L(M)$ of~\eqref{eq:lift}, the column indexed by $p$ is $t$.
Hence a vertex $i$ is adjacent to $p$ in $G^{\dagger}$ exactly when $t_i=1$, so the chords crossing $e$ are precisely the chords with $t_i=1$.
Before any reversal, the intersection matrix of the remaining chords is $K$, the principal block $L(M)[V]$.

Delete $e$ and reverse the order on $S_2$.
By Lemma~\ref{lem:reversal}, the crossing relation of two chords $i,j$ changes exactly when each of them has one endpoint on each of $S_1,S_2$, that is, exactly when $t_i=t_j=1$.
Consequently the off-diagonal entry of the reversed diagram is
\[
K_{ij}+t_it_j=M_{ij},
\]
and when $t_it_j=0$ the entry is unchanged, which is again $M_{ij}$.

Give sign $1$ exactly to the chords with endpoints on different arcs, that is, to the chords with $t_i=1$.
The resulting bouquet $B$ has a loop at $i$ precisely when $t_i=1$; since $t_i^2=t_i$ in $\F$, its diagonal entries are $M_{ii}$.
Therefore $M(I(B))=K+tt^T=M$, so $I(B)\cong G$.
The pair $(S_1,S_2)$ is a certificate: a chord with $t_i=0$ has both endpoints in the same arc, while a chord with $t_i=1$ meets both arcs.

Conversely, suppose $B$ has certificate $(S_1,S_2)$ and $I(B)\cong G$.
Let $\widehat B$ be the adjustment of Definition~\ref{def:pseudo}: reverse $S_2$, toggle the sign of every chord meeting both arcs, and add the edge $e=PQ$ with sign $0$.
By the certificate condition, an old chord meets both arcs exactly when it is twisted.
The toggling therefore makes every old sign zero, and the new edge also has sign zero, so $\widehat B$ is orientable.

For the crossings among old chords, Lemma~\ref{lem:reversal} shows that reversing $S_2$ changes the relation of chords $i,j$ exactly when both meet both arcs, that is, exactly when $t_i=t_j=1$.
Starting from the matrix $M=M(I(B))$, the new intersection matrix on the old chords is therefore $M+tt^T$.
The new chord $e$ crosses exactly the chords having one endpoint on each arc, which are exactly the originally twisted chords; hence its adjacency row is $t$, and $e$ carries no loop.
Thus the intersection matrix of $\widehat B$ is
\[
\begin{pmatrix}M+tt^T&t\\t^T&0\end{pmatrix}=L(M)=M(G^{\dagger}).
\]
Since $\widehat B$ is orientable, its signed chord diagram has all signs zero and is a circle representation of $I(\widehat B)$.
Hence $G^{\dagger}$ is a circle graph, so $G$ is stable.
\end{proof}

\begin{remark}\label{rem:existential}
The theorem asserts existence, not pseudo-orientability of every bouquet with intersection graph $G$.
Example~3.15 of~\cite{DK2026} gives two bouquets with identical binary intersection matrices and different pseudo-orientability properties.
Thus no invariance of pseudo-orientability under arbitrary changes of chord representation is used here.
\end{remark}

\section{Excluded toggle-minors}\label{sec:obstructions}
\subsection{A finite candidate family}
Fix a minor relation $\preceq$ on a class of graphs and a graph class $\mathcal K$.
The class $\mathcal K$ is \emph{minor-closed} if $H\in\mathcal K$ and $J\preceq H$ imply $J\in\mathcal K$.
An \emph{excluded minor class} for $\mathcal K$ is an equivalence class $\mathcal E$ contained in the complement of $\mathcal K$ such that every proper minor of every member of $\mathcal E$ lies in $\mathcal K$.
Here the relevant equivalence is toggle equivalence and ``proper'' means a strict decrease in order.
An \emph{excluded-minor representative} is one chosen graph from such a class; representatives need not be unique as labeled graphs.

\begin{figure}[htbp]
\centering
\begin{tikzpicture}[every node/.style={font=\scriptsize}]
\begin{scope}[scale=.62]
\foreach \a in {90,162,234,306,18}{\node[fgdot] (p\a) at (\a:1){};}
\draw[fgedge] (p90)--(p162)--(p234)--(p306)--(p18)--(p90);
\node[fgdot] (c) at (0,0){};
\foreach \a in {90,162,234,306,18}{\draw[fgedge] (c)--(p\a);}
\node at (0,-1.55) {$F_1$};
\end{scope}
\begin{scope}[shift={(3.05cm,0)},scale=.62]
\foreach \a in {90,30,-30,-90,-150,150}{\node[fgdot] (p\a) at (\a:1){};}
\draw[fgedge] (p90)--(p30)--(p-30)--(p-90)--(p-150)--(p150)--(p90);
\node[fgdot] (c) at (0,0){};
\foreach \a in {90,30,-30,-90,-150,150}{\draw[fgedge] (c)--(p\a);}
\node at (0,-1.55) {$F_2$};
\end{scope}
\begin{scope}[shift={(6.10cm,0)},scale=.62]
\foreach \k in {0,...,6}{\node[fgdot] (q\k) at ({90+\k*51.4286}:1){};}
\foreach \k in {0,...,6}{\pgfmathtruncatemacro{\n}{int(mod(\k+1,7))}\draw[fgedge] (q\k)--(q\n);}
\node[fgdot] (c) at (0,0){};
\foreach \k in {0,...,6}{\draw[fgedge] (c)--(q\k);}
\node at (0,-1.55) {$F_3$};
\end{scope}
\begin{scope}[shift={(9.15cm,0)},scale=.62]
\foreach \a in {90,30,-30,-90,-150,150}{\node[fgdot] (p\a) at (\a:1){};}
\draw[fgedge] (p90)--(p30)--(p-30)--(p-90)--(p-150)--(p150)--(p90);
\node[fgdot] (c) at (0,0){};
\draw[fgedge] (c)--(p90); \draw[fgedge] (c)--(p-30); \draw[fgedge] (c)--(p-150);
\node at (0,-1.55) {$F_4$};
\end{scope}
\begin{scope}[shift={(12.20cm,0)},scale=.62]
\node[fgdot] (p1) at (0,1){}; \node[fgdot] (p3) at (0.95,-0.45){};
\node[fgdot] (p5) at (-0.95,-0.45){};
\node[fgdot] (p2) at (0.72,0.5){}; \node[fgdot] (p4) at (0,-0.68){};
\node[fgdot] (p6) at (-0.72,0.5){}; \node[fgdot] (p7) at (0,0.02){};
\draw[fgedge] (p1)--(p3); \draw[fgedge] (p3)--(p4)--(p5); \draw[fgedge] (p1)--(p5);
\draw[fgedge] (p1)--(p2)--(p3); \draw[fgedge] (p1)--(p6)--(p5);
\draw[fgedge] (p7)--(p1); \draw[fgedge] (p7)--(p3); \draw[fgedge] (p7)--(p5);
\node at (0,-1.55) {$F_5$};
\end{scope}
\begin{scope}[shift={(0,-3.05cm)},scale=.62]
\node[fgdot] (p1) at (0,1){}; \node[fgdot] (p3) at (0.87,-0.5){};
\node[fgdot] (p5) at (-0.87,-0.5){};
\draw[fgedge] (p1)--(p3)--(p5)--(p1);
\node[fgdot] (p2) at (0.38,0.55){}; \node[fgdot] (p4) at (0,-0.25){};
\node[fgdot] (p6) at (-0.38,0.55){}; \node[fgdot] (p7) at (0,0.12){};
\draw[fgedge] (p2)--(p1); \draw[fgedge] (p2)--(p3);
\draw[fgedge] (p4)--(p3); \draw[fgedge] (p4)--(p5);
\draw[fgedge] (p6)--(p1); \draw[fgedge] (p6)--(p5);
\draw[fgedge] (p7)--(p2); \draw[fgedge] (p7)--(p4); \draw[fgedge] (p7)--(p6);
\node at (0,-1.55) {$F_6$};
\end{scope}
\begin{scope}[shift={(3.05cm,-3.05cm)},scale=.62]
\foreach \a in {90,30,-30,-90,-150,150}{\node[fgdot] (p\a) at (\a:1){};}
\draw[fgedge] (p90)--(p30)--(p-30)--(p-90)--(p-150)--(p150)--(p90);
\node[fgdot] (c) at (0,0){};
\foreach \a in {90,30,-90,-150,150}{\draw[fgedge] (c)--(p\a);}
\node at (0,-1.55) {$F_7$};
\end{scope}
\begin{scope}[shift={(6.10cm,-3.05cm)},scale=.62]
\foreach \x in {-1,0,1}{\foreach \y in {-1,0,1}{\node[fgdot] (g\x\y) at (\x,\y){};}}
\foreach \x in {-1,0,1}{\draw[fgedge] (g\x-1)--(g\x0)--(g\x1);}
\foreach \y in {-1,0,1}{\draw[fgedge] (g-1\y)--(g0\y)--(g1\y);}
\node at (0,-1.55) {$F_8$};
\end{scope}
\begin{scope}[shift={(9.15cm,-3.05cm)},scale=.62]
\node[fgdot] (a) at (0,1.05){}; \node[fgdot] (b) at (0.95,-0.55){};
\node[fgdot] (c) at (-0.95,-0.55){}; \node[fgdot] (d) at (0,0.10){};
\node[fgdot] (ab) at (0.475,0.25){}; \node[fgdot] (ac) at (-0.475,0.25){};
\node[fgdot] (bc) at (0,-0.55){}; \node[fgdot] (cd) at (-0.475,-0.225){};
\node[fgdot] (bd) at (0.475,-0.225){}; \node[fgdot] (ad) at (0,0.575){};
\draw[fgedge] (a)--(ab)--(b); \draw[fgedge] (b)--(bc)--(c); \draw[fgedge] (c)--(ac)--(a);
\draw[fgedge] (c)--(cd)--(d); \draw[fgedge] (b)--(bd)--(d); \draw[fgedge] (a)--(ad)--(d);
\node at (0,-1.55) {$F_9$};
\end{scope}
\begin{scope}[shift={(12.20cm,-3.05cm)},scale=.62]
	\node[fgdot] (p1) at (-0.70,0.50){};
	\node[fgdot] (p4) at (-1.20,-0.35){};
	\node[fgdot] (p5) at (-0.18,-0.35){};
	\draw[fgedge] (p1)--(p4)--(p5)--(p1);
	\node[fgdot] (p2) at (0.70,0.50){};
	\node[fgdot] (p3) at (1.20,-0.35){};
	\node[fgdot] (p7) at (0.18,-0.35){};
	\draw[fgedge] (p2)--(p3)--(p7)--(p2);
	\draw[fgedge] (p1)--(p2);
	\node[fgdot] (p6) at (0,-1.05){};
	\draw[fgedge] (p5)--(p6)--(p7);
	\draw[fgedge] (p3) to[bend left=14] (p4);
	\node at (0,-1.55) {$F_{10}$};
\end{scope}
\begin{scope}[shift={(0,-6.10cm)},scale=.62]
\foreach \a in {90,30,-30,-90,-150,150}{\node[fgdot] (p\a) at (\a:1){};}
\draw[fgedge] (p90)--(p30)--(p-30)--(p-90)--(p-150)--(p150)--(p90);
\node[fgdot] (p7) at (-0.35,0.02){}; \node[fgdot] (p8) at (0.35,-0.45){};
\draw[fgedge] (p7)--(p90); \draw[fgedge] (p7)--(p-150);
\draw[fgedge] (p8)--(p-30); \draw[fgedge] (p8)--(p-90); \draw[fgedge] (p8)--(p7);
\node at (0,-1.55) {$F_{11}$};
\end{scope}
\begin{scope}[shift={(3.05cm,-6.10cm)},scale=.62]
\foreach \k in {0,...,7}{\node[fgdot] (p\k) at ({90+\k*45}:1){};}
\foreach \k in {0,...,7}{\pgfmathtruncatemacro{\n}{int(mod(\k+1,8))}\draw[fgedge] (p\k)--(p\n);}
\node[fgdot] (c) at (0,0){};
\foreach \k in {0,1,3,4,5,7}{\draw[fgedge] (c)--(p\k);}
\draw[fgedge] (p0) to[bend right=14] (p4);
\node at (0,-1.55) {$F_{12}$};
\end{scope}
\begin{scope}[shift={(6.10cm,-6.10cm)},scale=.62]
\foreach \k in {0,...,6}{\node[fgdot] (p\k) at ({90+\k*51.4286}:1){};}
\foreach \k in {0,...,6}{\pgfmathtruncatemacro{\n}{int(mod(\k+1,7))}\draw[fgedge] (p\k)--(p\n);}
\node[fgdot] (c) at (0,0){};
\foreach \k in {0,2,5}{\draw[fgedge] (c)--(p\k);}
\node at (0,-1.55) {$F_{13}$};
\end{scope}
\begin{scope}[shift={(9.15cm,-6.10cm)},scale=.62]
\foreach \a in {90,30,-30,-90,-150,150}{\node[fgdot] (p\a) at (\a:1){};}
\draw[fgedge] (p90)--(p30)--(p-30)--(p-90)--(p-150)--(p150)--(p90);
\node[fgdot] (c) at (0,0){};
\foreach \a in {-30,-90,-150,150}{\draw[fgedge] (c)--(p\a);}
\node[fgdot] (s) at (-0.5,-0.55){};
\draw[fgedge] (s)--(p-90); \draw[fgedge] (s)--(p-150); \draw[fgedge] (s)--(c);
\node at (0,-1.55) {$F_{14}$};
\end{scope}
\begin{scope}[shift={(12.20cm,-6.10cm)},scale=.62]
\foreach \k in {0,...,7}{\node[fgdot] (p\k) at ({90+\k*45}:1){};}
\foreach \k in {0,...,7}{\pgfmathtruncatemacro{\n}{int(mod(\k+1,8))}\draw[fgedge] (p\k)--(p\n);}
\draw[fgedge] (p0) to[bend left=14] (p4);
\draw[fgedge] (p1) to[bend left=14] (p5);
\draw[fgedge] (p2) to[bend left=14] (p6);
\draw[fgedge] (p3) to[bend left=14] (p7);
\node at (0,-1.55) {$F_{15}$};
\end{scope}
\end{tikzpicture}
\caption{The fifteen excluded pivot-minors $\cF$ for circle graphs, redrawn after Figure~2 of~\cite{GO2009}. The subscript $i$ follows the numbering of the input matrices in the ancillary file \texttt{anc/enumeration/inputs.json}, which is the row-major order of that figure. Every $F_i$ is connected and is not a circle graph.}
\label{fig:excluded}
\end{figure}
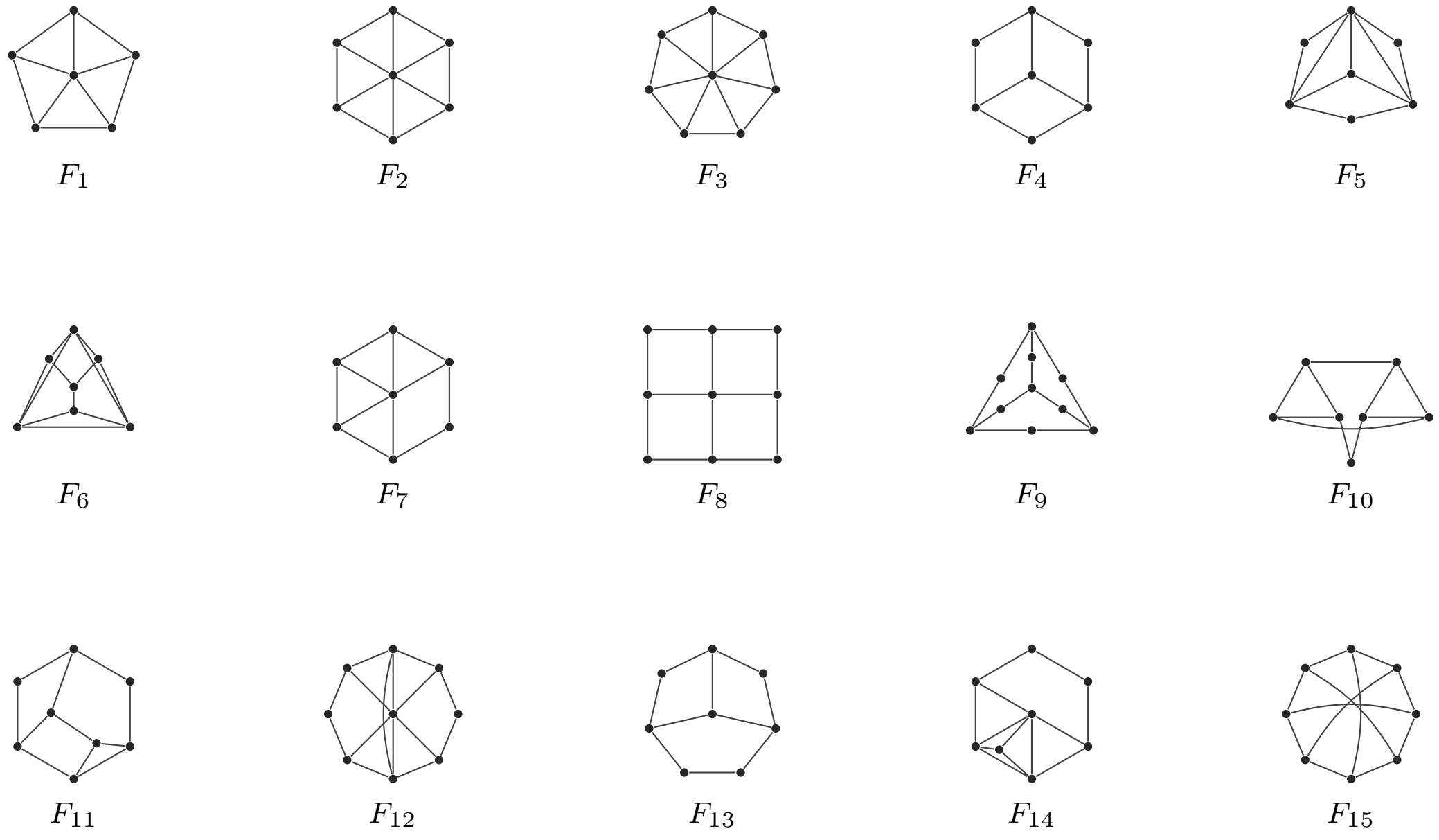

Fix one representative from each of the fifteen excluded pivot-minor classes for circle graphs, as in~\cite[Theorem 1.2 and Figure 2]{GO2009}, and call the resulting family $\cF$.
Its defining property is
\begin{equation}\label{eq:GO}
H\text{ is a circle graph}\quad\Longleftrightarrow\quad
\text{no }F\in\cF\text{ is a pivot-minor of }H.
\end{equation}
Every member of $\cF$ is connected and is not a circle graph.
The labeled adjacency matrices used in our enumeration are specified in the ancillary file \texttt{anc/enumeration/inputs.json}.
The fifteen graphs are drawn in Figure~\ref{fig:excluded}.

\begin{definition}[Candidates and minimal classes]\label{def:candidates}
For $F\in\cF$ and a nonempty set $N\subseteq V(F)$, let $F_N^+$ be the simple graph obtained by adjoining a new vertex $r$ with neighborhood exactly $N$.
Form a finite family $\cC$ by taking the following graphs and retaining only those that are looped circle graphs:
\begin{enumerate}[label=\textup{(\arabic*)}]
\item $\Rroot Fr$ for every $F\in\cF$ and $r\in V(F)$;
\item $\Rroot{F_N^+}r$ for every $F\in\cF$ and every nonempty $N\subseteq V(F)$.
\end{enumerate}
Identify candidates that are toggle-equivalent.
Among the resulting classes retain the minimal ones under $\preceq$, and choose one representative from each. Their family is denoted by $\cG$.
\end{definition}

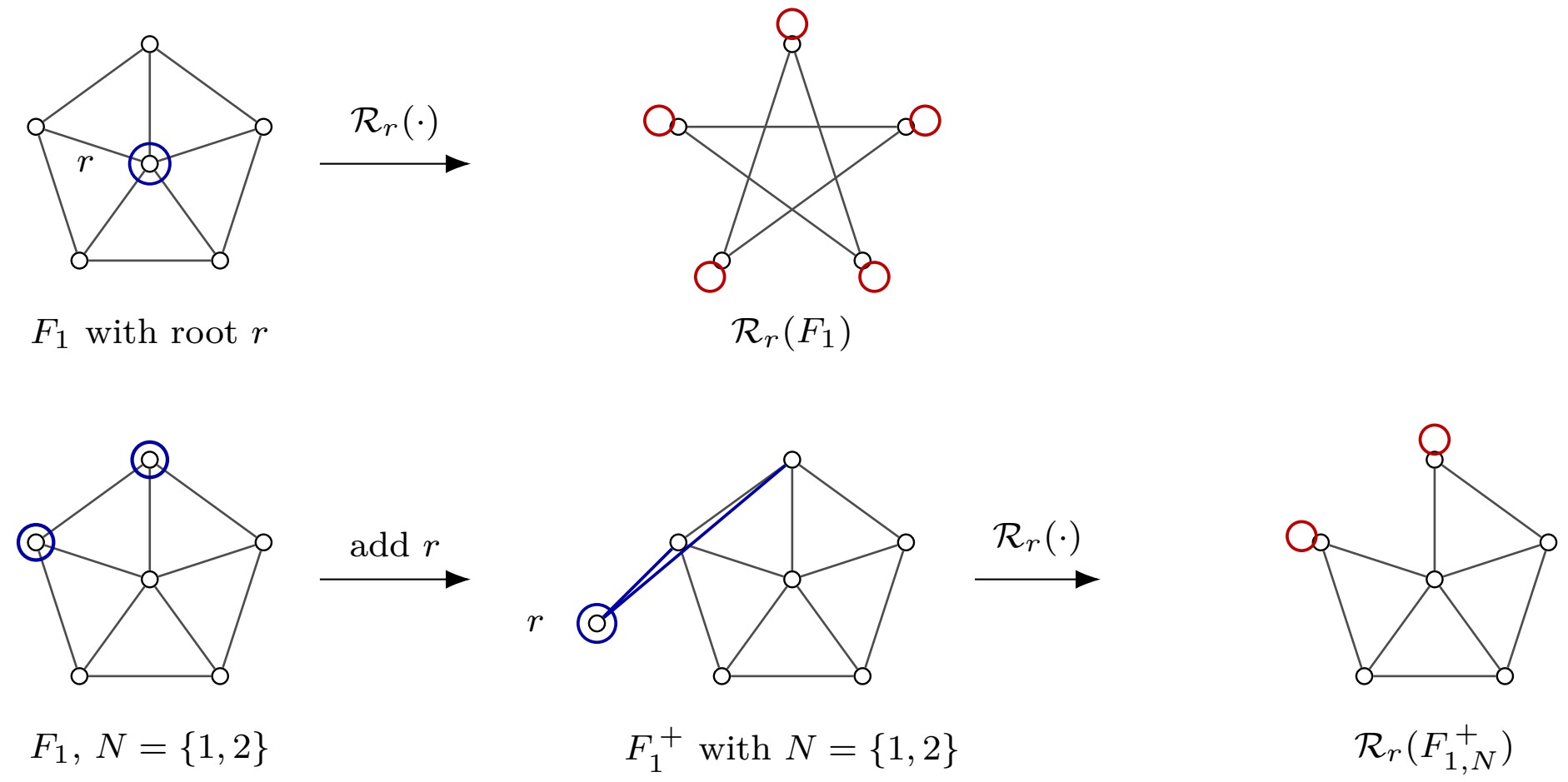
\begin{figure}[htbp]
\centering
\begin{tikzpicture}[scale=1.0,every node/.style={font=\scriptsize}]
\begin{scope}
\foreach \a in {90,162,234,306,18}{\node[gv] (a\a) at (\a:.95){};}
\draw[ge] (a90)--(a162)--(a234)--(a306)--(a18)--(a90);
\node[gv] (ar) at (0,0){};
\foreach \a in {90,162,234,306,18}{\draw[ge] (ar)--(a\a);}
\draw[blue!65!black,line width=.7pt] (ar) circle[radius=.16];
\node[left=0.3cm] at (ar) {$r$};
\node at (0,-1.35) {$F_1$ with root $r$};
\end{scope}
\draw[-{Latex[length=2mm]}] (1.35,0)--(2.55,0);
\node[above] at (1.95,0.05) {$\mathcal R_{r}(\cdot)$};
\begin{scope}[xshift=5.1cm]
\foreach \a in {90,162,234,306,18}{\node[gv] (b\a) at (\a:.95){};}
\foreach \a/\b in {90/234,234/18,18/162,162/306,306/90}{\draw[ge] (b\a)--(b\b);}
\foreach \a in {90,162,234,306,18}{\draw[le] ($(b\a)+(\a:.16)$) circle[radius=.115];}
\node at (0,-1.35) {$\mathcal R_{r}(F_1)$};
\end{scope}
\begin{scope}[yshift=-3.3cm]
\foreach \a in {90,162,234,306,18}{\node[gv] (c\a) at (\a:.95){};}
\draw[ge] (c90)--(c162)--(c234)--(c306)--(c18)--(c90);
\node[gv] (cr) at (0,0){};
\foreach \a in {90,162,234,306,18}{\draw[ge] (cr)--(c\a);}
\draw[blue!65!black,line width=.8pt] (c90) circle[radius=.14];
\draw[blue!65!black,line width=.8pt] (c162) circle[radius=.14];
\node at (0,-1.35) {$F_1$, $N=\{1,2\}$};
\end{scope}
\draw[-{Latex[length=2mm]}] (1.35,-3.3)--(2.55,-3.3);
\node[above] at (1.95,-3.25) {add $r$};
\begin{scope}[xshift=5.1cm,yshift=-3.3cm]
\foreach \a in {90,162,234,306,18}{\node[gv] (d\a) at (\a:.95){};}
\draw[ge] (d90)--(d162)--(d234)--(d306)--(d18)--(d90);
\node[gv] (dr) at (0,0){};
\foreach \a in {90,162,234,306,18}{\draw[ge] (dr)--(d\a);}
\node[gv] (dn) at (-1.55,-0.35){};
\draw[ne] (dn)--(d90); \draw[ne] (dn)--(d162);
\draw[blue!65!black,line width=.7pt] (dn) circle[radius=.15];
\node[left=0.28cm] at (dn) {$r$};
\node at (0,-1.35) {$F_1^{\,+}$ with $N=\{1,2\}$};
\end{scope}
\draw[-{Latex[length=2mm]}] (6.55,-3.3)--(7.55,-3.3);
\node[above] at (7.05,-3.25) {$\mathcal R_{r}(\cdot)$};
\begin{scope}[xshift=10.2cm,yshift=-3.3cm]
\foreach \a in {90,162,234,306,18}{\node[gv] (e\a) at (\a:.95){};}
\draw[ge] (e162)--(e234)--(e306)--(e18)--(e90);
\node[gv] (er) at (0,0){};
\foreach \a in {90,162,234,306,18}{\draw[ge] (er)--(e\a);}
\draw[le] ($(e90)+(90:.16)$) circle[radius=.115];
\draw[le] ($(e162)+(162:.16)$) circle[radius=.115];
\node at (0,-1.35) {$\mathcal R_{r}(F_{1,N}^{\,+})$};
\end{scope}
\end{tikzpicture}
\caption{The two rooted constructions of Definition~\ref{def:candidates} on $F_1$. (a) Type~(1): root $F_1$ at its centre $r$ and apply the rooted inverse. Since $t=\mathbf 1$, the result is $K+tt^T$ with $K$ the adjacency matrix of the rim $5$-cycle, namely the all-looped complement of that cycle. (b) Type~(2) with $N=\{1,2\}$: adjoin a new vertex $r$ adjacent exactly to $N$ (blue), then apply the rooted inverse; the result is $F_1$ with the edge $12$ toggled off and loops added at $1$ and $2$. Both displayed candidates have a circle graph as underlying graph, so both belong to $\cC$.}
\label{fig:candidates}
\end{figure}

Figure~\ref{fig:candidates} illustrates the two constructions on $F_1$.

If $m=\max_{F\in\cF}|V(F)|$, all candidates have order at most $m$.
Every candidate is non-stable.
Indeed, for a type~(1) candidate $\Rroot Fr$, Proposition~\ref{prop:inverse} gives
$(\Rroot Fr)^{\dagger}\cong F$ up to relabeling of the new vertex, so its lift is isomorphic to $F$.
For a type~(2) candidate $\Rroot{F_N^+}r$, the same proposition gives
$(\Rroot{F_N^+}r)^{\dagger}\cong F_N^+$, and $F_N^+$ contains $F$ as the induced subgraph on $V(F)$.
Since $F$ is not a circle graph and circle graphs are closed under induced subgraphs, neither lift is a circle graph.
Thus no candidate is stable.

\begin{lemma}[Completeness of the candidate construction]\label{lem:complete}
Every non-stable looped circle graph has a toggle-minor isomorphic to a member of $\cC$.
\end{lemma}
\begin{proof}
Let $G$ be non-stable. By~\eqref{eq:GO}, there is a simple graph $H$ pivot-equivalent to $G^{\dagger(p)}$ and a set $S$ with $H[S]\cong F\in\cF$.
Transport the root label and put $G'=\Rroot Hp$.
By Corollary~\ref{cor:projection}, $G'$ is toggle-equivalent to $G$, so it is a looped circle graph.

Write $M(H)=\left(\begin{smallmatrix}K&t\\t^T&0\end{smallmatrix}\right)$, with the root $p$ last.
Then $G'$ has matrix $K+tt^T$ on $V(H)\setminus\{p\}$.

Suppose first that $p\in S$.
Restricting the matrix of $G'$ to $S\setminus\{p\}$ gives
\[
(K+tt^T)[S\setminus\{p\}]
=K[S\setminus\{p\}]+t[S\setminus\{p\}]\,t[S\setminus\{p\}]^T,
\]
which is exactly the matrix of the rooted inverse $\Rroot{H[S]}p$.
Thus
\[
G'[S\setminus\{p\}]=\Rroot{H[S]}p.
\]
Since $H[S]\cong F$, this is isomorphic to a type~(1) candidate $\Rroot Fr$.
It passes the looped-circle-graph filter because it is an induced subgraph of $G'$.

Now suppose that $p\notin S$.
If $N_H(p)\cap S=\varnothing$, then $t[S]=0$, and~\eqref{eq:root} gives
\[
G'[S]=K[S]+t[S]t[S]^T=K[S]=H[S]\cong F.
\]
This is impossible: $G'[S]$ is an induced subgraph of the looped circle graph $G'$, so its underlying simple graph is a circle graph, whereas $F$ is not.

Hence $N=N_H(p)\cap S$ is nonempty.
The induced subgraph $H[S\cup\{p\}]$ is isomorphic to $F_N^+$: its part on $S$ is $H[S]\cong F$, and the new vertex $p$ has neighborhood exactly $N$ in $S$.
Therefore
\[
G'[S]=K[S]+t[S]t[S]^T=\Rroot{H[S\cup\{p\}]}p
\]
is isomorphic to a type~(2) candidate $\Rroot{F_N^+}r$, and it again passes the filter because it is an induced subgraph of $G'$.

In both cases $G'[S]$ or $G'[S\setminus\{p\}]$ is an induced subgraph of $G'$ and hence a toggle-minor of $G'$; by transitivity of $\preceq$ and the toggle equivalence of $G'$ and $G$, it is a toggle-minor of $G$.
\end{proof}

\begin{lemma}[Finite descent]\label{lem:descent}
Every candidate in $\cC$ has a toggle-minor isomorphic to a member of $\cG$.
\end{lemma}
\begin{proof}
If its class is not minimal among candidate classes, pass to a strictly smaller candidate class below it.
Distinct comparable classes of the same order cannot occur, because an equal-order toggle-minor involves no vertex deletion and hence is toggle-equivalent to its parent.
Every strict step decreases the order, which is a nonnegative integer, so a strictly descending chain of candidate classes is finite and the process terminates at a minimal candidate class.
Transitivity of $\preceq$ then shows that the terminal class is a toggle-minor of the original candidate, which proves the assertion.
\end{proof}

\begin{theorem}[Excluded toggle-minor characterization]\label{thm:forbidden}
Let $G$ be a looped circle graph and let $\cG$ be as in Definition~\ref{def:candidates}.
Then $G$ is stable if and only if it has no toggle-minor isomorphic to a member of $\cG$.
Moreover, $\cG$ is a complete set of representatives of the minimal non-stable toggle-equivalence classes within the looped circle graphs.
\end{theorem}
\begin{proof}
If $G$ is stable, all its toggle-minors are stable by Proposition~\ref{prop:closure}, whereas no member of $\cG$ is stable.
If $G$ is not stable, Lemma~\ref{lem:complete} gives a candidate minor and Lemma~\ref{lem:descent} gives a member of $\cG$ below that candidate.

For minimality, suppose $J\in\cG$ has a non-stable proper toggle-minor $J'$.
It is a looped circle graph, so Lemma~\ref{lem:complete} gives a candidate $C\preceq J'$.
Then $C\preceq J$ and $|V(C)|\le |V(J')|<|V(J)|$, contradicting the minimality of $J$ among candidate classes.
Conversely, any minimal non-stable class contains a member of $\cG$ as a minor by the first part.
Its minimality forces equal order and therefore toggle equivalence.
\end{proof}

\subsection{Exact recognition and computational certificates}\label{sec:computation}
We give precise meanings to the computational tests and certificates used below.
The \emph{isomorphism canonical form} of a binary symmetric matrix is a chosen representative of its orbit under simultaneous row-and-column permutations, chosen so that two matrices have the same form exactly when they are isomorphic.
The implementation uses an exhaustive vertex-ordering search, retaining the diagonal.
Its partition refinement groups vertices according to their loop status and their numbers of neighbors in each current block, repeatedly splitting blocks when these counts differ.
Individualization means separating one vertex from a nonsingleton block and exploring all resulting choices.
Branches corresponding to twins may be identified: two vertices are twins here if they have the same loop status and the same adjacency to every vertex outside the pair; their transposition is then an automorphism, that is, an isomorphism from the graph to itself.
The \emph{pivot-isomorphism orbit} of $M$ is
\[
\mathcal O(M)=\{\operatorname{can}(M*X):X\subseteq V,\ \det M[X]=1\}.
\]
Here $\operatorname{can}$ denotes the canonical form, and the empty pivot is included.
The pivot composition law ensures that no additional states arise from longer transform sequences.

\begin{definition}[The Naji system]\label{def:naji}
For a simple graph $H$, introduce variables $\beta(v,w)\in\F$ for all ordered pairs of distinct vertices.
Its \emph{Naji system} consists of
\begin{align*}
\beta(v,w)+\beta(w,v)&=1 &&(vw\in E(H)),\\
\beta(x,v)+\beta(x,w)&=0 &&(vw\in E(H),\ xv,xw\notin E(H)),\\
\beta(v,w)+\beta(w,v)+\beta(x,v)+\beta(x,w)&=1
&& (xv,xw\in E(H),\ vw\notin E(H)),
\end{align*}
where the vertices in each triple are distinct.
\end{definition}
Naji's theorem states that this system is consistent if and only if $H$ is a circle graph~\cite[Theorem 1.1]{GL2020}.
For a looped graph the test is applied after removing the diagonal of its adjacency matrix.

A \emph{positive linear-system certificate} is a binary assignment to all variables satisfying every equation.
A \emph{negative certificate} is a set of equations whose sum over $\F$ has zero left-hand side and right-hand side $1$.
For systems $Ax=b$ over a field, Gaussian elimination produces one of these two certificates.
Both are checked by direct substitution or addition, without rerunning the original search.

An \emph{equivalence certificate} records a feasible pivot set and a simultaneous row-and-column permutation.
A \emph{minor certificate} records a feasible pivot, an ordered list of retained indices, and an isomorphism; the target may also be specified by an equivalent pivoted representative.
Checking these data consists of evaluating~\eqref{eq:ppt}, taking the indicated principal submatrix, and comparing the resulting labeled matrices.

The exact enumeration uses the following finite procedure:
\begin{enumerate}
\item Generate both types in Definition~\ref{def:candidates} from the fifteen input matrices.
\item Remove duplicate isomorphism types and solve each underlying simple graph's Naji system by exact binary elimination.
\item Enumerate every feasible pivot, including the empty one, and merge toggle-equivalent candidates using their canonical forms.
\item Remove each candidate class containing a strictly lower-order candidate class; save a minor certificate for each removal.
\end{enumerate}
No timeout is interpreted as an inconsistency or a noncontainment result.

\begin{algorithm}[htbp]
\caption{Candidate generation, toggle classification, and minimal representatives}\label{alg:representatives}
\begin{algorithmic}[1]
\Require The family $\cF$ of fifteen excluded circle-graph pivot-minors
\Ensure One representative per minimal non-stable toggle class
\State $\cC\gets\varnothing$
\For{$F\in\cF$}
  \For{$r\in V(F)$}
    \State $A\gets\Rroot Fr$
    \If{$A$ is a looped circle graph} \State $\cC\gets\cC\cup\{A\}$ \EndIf
  \EndFor
  \For{$\varnothing\ne N\subseteq V(F)$}
    \State Construct $F_N^+$ with new vertex $r$; set $A\gets\Rroot{F_N^+}r$
    \If{$A$ is a looped circle graph} \State $\cC\gets\cC\cup\{A\}$ \EndIf
  \EndFor
\EndFor
\State Remove isomorphic duplicates from $\cC$, preserving loops
\State Enumerate all feasible pivots, including $\varnothing$, and merge candidates
\Statex\hspace{\algorithmicindent} into the family $\mathcal T$ of toggle-equivalence classes
\State $\cG\gets\varnothing$
\For{$\mathcal A\in\mathcal T$ in nondecreasing order of graph order}
  \State Choose a representative $A\in\mathcal A$
  \If{no $J\in\cG$ is a proper toggle-minor of $A$}
    \State $\cG\gets\cG\cup\{A\}$
  \Else
    \State Save a pivot, retained indices, and isomorphism witnessing $J\preceq A$
  \EndIf
\EndFor
\State \Return $\cG$
\end{algorithmic}
\end{algorithm}

Algorithm~\ref{alg:representatives} combines both candidate constructions with global minimization.
The recognition and minor tests are exact; the latter enumerate feasible pivots and induced subsets of the required order.
Its correctness follows from Lemmas~\ref{lem:complete} and~\ref{lem:descent} and Theorem~\ref{thm:forbidden}.
Processing by order ensures that a discarded class already contains a retained smaller minimal class, so it need not be reconsidered.
All loops range over finite sets.

\begin{proposition}[Exact enumeration]\label{prop:enumeration}
The family $\cG$ has $22$ members: one of order $5$, twelve of order $6$, and nine of order $7$.
One choice of representatives is displayed in Appendix~\ref{app:representatives}.
\end{proposition}
\begin{proof}[Computational proof and verification]
The input is checked against Figure~2 of~\cite{GO2009}.
Each of the fifteen inputs is certified not to be a circle graph; every one-vertex deletion from every graph in its pivot-isomorphism orbit is certified to be a circle graph.
The inputs belong to pairwise distinct pivot-equivalence classes.
Together with the fifteen-obstruction theorem, these checks ensure coverage of the input obstruction classes.
The exact stage counts are
\begin{center}
\begin{tabular}{lr}\toprule
Stage&Count\\\midrule
Rooted constructions, before deduplication&4261\\
Distinct isomorphism types&1338\\
Looped circle graphs retained&698\\
Toggle-equivalence classes&237\\
Nonminimal classes removed&215\\
Minimal classes retained&22\\\bottomrule
\end{tabular}
\end{center}
All recognition, equivalence, and removal certificates are saved.
An additional direct test certifies that each retained representative has a non-circle-graph lift, while every one-vertex deletion from every member of its toggle orbit has a circle-graph lift.
There are $2589$ such deletion tests, all affirmative.
Since stability is minor-closed, these checks establish minimality independently of the particular removal certificates.
The equal-order orbits of different retained representatives are disjoint.

For a separate check on the recognizer, all $10395$ pairings of twelve cyclically ordered endpoints are enumerated.
They give $154$ isomorphism types of six-vertex circle graphs; comparison against all $32768$ labeled six-vertex simple graphs agrees with the Naji test in every case.
Finally, an actual double-occurrence word is generated for each of the $22$ representatives and its crossing relations are checked pair by pair.
The candidate completeness and descent arguments in Theorem~\ref{thm:forbidden} turn this finite calculation into the claimed complete list.
\end{proof}

\begin{remark}[Data and reproducibility]
The ancillary directory \texttt{anc/enumeration/} contains the fifteen inputs, the exact programs, representative matrices, double-occurrence words, certificates, and verification summaries.
Its \texttt{README.md} describes the files and gives commands for reproducing and checking the computation.
The original transcription of the sixth input had an incorrect central neighborhood; the input file records the correction from the published figure.
The program generates all candidates directly from the fifteen input matrices.
In the JSON files, vertices are indexed by $0,\ldots,n-1$; in the appendix drawings, the corresponding vertices are labeled $1,\ldots,n$.
\end{remark}

\section{Interlace polynomials and quasi-tree counting}\label{sec:polynomial}
The determinant identity of Lemma~\ref{lem:parity} was the algebraic engine of Sections~\ref{sec:lift}--\ref{sec:obstructions}: it compares the principal minors of the lift $L(M)$ with those of $M$.
This section proves its nullity counterpart (Lemma~\ref{lem:nullity}) and derives the resulting polynomial and enumerative statements.
The identity of Theorem~\ref{thm:polynomial} is stated for every looped simple graph, since no circle-graph hypothesis is needed; stability enters only in Corollary~\ref{cor:determinant}.
Apart from the use of Theorem~\ref{thm:realization} in Corollary~\ref{cor:quasitree}, the section does not depend on Section~\ref{sec:obstructions} and may be read immediately after Section~\ref{sec:bouquets}.
The counting results below concern the spanning quasi-trees introduced in Definition~\ref{def:ribbon}.
\subsection{A nullity pair}
Two symmetric matrices over a field are \emph{congruent} if one is $P^TAP$ for the other and some invertible matrix $P$; congruence preserves rank.
The \emph{direct sum} $A\oplus B$ is the block-diagonal matrix with blocks $A,B$.
The bilinear form of a symmetric matrix $K$ is $(x,y)\mapsto x^TKy$.
Its \emph{radical} is $\{u:u^TKv=0\text{ for all }v\}$, which equals $\ker K$ because $K$ is symmetric.
For a subspace $U$, write $U^{\perp}=\{v:u^Tv=0\text{ for all }u\in U\}$ for its orthogonal complement with respect to the standard dot product.
A subspace $W$ with $V=W\oplus\ker K$ whose restricted form is nondegenerate, meaning that its radical in $W$ is zero, is called a \emph{nondegenerate part} of $K$.
Such a part exists: any complement of $\ker K$ maps isomorphically to $V/\ker K$, and the form induced on $V/\ker K$ is nondegenerate because the radical of $K$ is $\ker K$.

\begin{lemma}\label{lem:nullity}
Let $M$ be symmetric over $\F$, $t=\diag M$, and $K=M+tt^T$.
As unordered pairs with multiplicity,
\begin{equation}\label{eq:nullity}
\{\nul K,\nul L(M)\}=\{\nul M,\nul M+1\}.
\end{equation}
\end{lemma}
\begin{proof}
If $t\in\im K$, choose $v$ with $Kv=t$.
Since $K$ is alternating, $t^Tv=v^TKv=0$.
Then
\[
M=K+tt^T=K+Kvt^T=K(I+vt^T).
\]
The matrix $I+vt^T$ is invertible: indeed
\[
(I+vt^T)^2=I+2vt^T+v(t^Tv)t^T=I
\]
in $\F$, so $(I+vt^T)^{-1}=I+vt^T$.
Hence $\rk M=\rk K$.

For the lift, put
\[
U=\begin{pmatrix}I&v\\0&1\end{pmatrix}.
\]
Using $Kv=t$ and $t^Tv=0$, we compute
\[
U^T L(M)U
=\begin{pmatrix}K&0\\0&0\end{pmatrix}.
\]
Thus $L(M)$ is congruent to $K\oplus[0]$, so
\[
\nul L(M)=\nul K+1=\nul M+1.
\]
This is the first case.

Now suppose $t\notin\im K$.
Since $K$ is symmetric, $\im K=(\ker K)^{\perp}$; hence there is $z\in\ker K$ with $t^Tz=1$.
Choose a basis $w_1,\ldots,w_r$ of a nondegenerate part of $K$, where $r=\rk K$, and replace each $w_i$ by
\[
w_i'=w_i+(t^Tw_i)z.
\]
Then $t^Tw_i'=t^Tw_i+(t^Tw_i)(t^Tz)=0$.
Moreover, because $Kz=0$,
\[
K(w_i',w_j')=K(w_i,w_j)
\]
for all $i,j$.
Since each $w_i'$ differs from $w_i$ by an element of $\ker K$, the $w_i'$ remain linearly independent modulo $\ker K$ and hence span a complement of $\ker K$; by definition, this complement is a nondegenerate part.
Complete this basis by a basis $z_1=z,z_2,\ldots,z_d$ of $\ker K$, chosen so that $t^Tz_1=1$ and $t^Tz_j=0$ for $j\ge2$; such a choice is possible because $t$ restricts to a nonzero linear functional on $\ker K$.
In this basis,
\[
K=K_0\oplus 0_d,\qquad t=(0,\ldots,0,1,0,\ldots,0)^T,
\]
where $K_0$ is invertible and the single $1$ of $t$ lies in the zero block.
Consequently $M=K+tt^T$ is block diagonal with blocks $K_0$ and $E_{11}$, the latter of rank $1$; hence
\[
\rk M=\rk K+1.
\]
For the lift, the same basis exhibits
\[
L(M)=\begin{pmatrix}
K_0&0&0\\
0&0_d&e_1\\
0&e_1^T&0
\end{pmatrix},
\]
where $e_1$ is the first standard basis vector of the zero block.
The block on the radical coordinates and the new coordinate is, after reordering, the direct sum of the hyperbolic block
$\left(\begin{smallmatrix}0&1\\1&0\end{smallmatrix}\right)$ on $z_1,p$ and a zero block of order $d-1$, so it has rank $2$.
Thus $\rk L(M)=\rk K+2$, and
\[
\nul L(M)=(n+1)-(\rk K+2)=\nul K-1.
\]
Since $\rk M=\rk K+1$, we also have $\nul M=n-(\rk K+1)=\nul K-1$.
Therefore $\nul K=\nul M+1$ and $\nul L(M)=\nul M$.

The empty matrix satisfies the formula directly: $\nul M=\nul K=0$ and $\nul L(M)=1$.
In both cases the unordered pairs coincide, as claimed.
\end{proof}

\begin{definition}[Vertex-nullity interlace polynomial]\label{def:q}
For a looped simple graph $G$ and an indeterminate $y$, define the integer-coefficient polynomial
\[
q(G;y)=\sum_{X\subseteq V(G)}(y-1)^{\nul M(G)[X]}.
\]
All nullities are computed over $\F$.
This is the vertex-nullity specialization of the interlace polynomial~\cite{ABS2004}.
\end{definition}

\begin{theorem}[Lift identity]\label{thm:polynomial}
For every looped simple graph $G$,
\begin{equation}\label{eq:q}
q(G^{\dagger(p)};y)=yq(G;y).
\end{equation}
\end{theorem}
\begin{proof}
The subsets of $V(G)\cup\{p\}$ partition into pairs $X$, $X\cup\{p\}$ indexed by $X\subseteq V(G)$.
For such an $X$, the principal matrix of $G^{\dagger}$ on $X$ is $K[X]$, and the principal matrix on $X\cup\{p\}$ is
\[
L(M)[X\cup\{p\}]=\begin{pmatrix}K[X]&t[X]\\t[X]^T&0\end{pmatrix}=L(M[X]).
\]
Put $d_X=\nul M[X]$.
By Lemma~\ref{lem:nullity} applied to the matrix $M[X]$, whose diagonal vector is $t[X]$ and whose associated matrix is $K[X]=M[X]+t[X]t[X]^T$, the two nullities are $d_X$ and $d_X+1$ in some order.
Their combined contribution to $q(G^{\dagger};y)$ is therefore
\[
(y-1)^{d_X}+(y-1)^{d_X+1}=y(y-1)^{d_X}.
\]
Summing over all $X\subseteq V(G)$ gives
\[
q(G^{\dagger};y)=\sum_{X\subseteq V(G)}y(y-1)^{d_X}=y\,q(G;y).
\]
This is an identity in $\mathbb Z[y]$, so it holds for every value of $y$, including $y=1$; at $y=1$ both sides equal the number of $X\subseteq V(G)$ with $d_X=0$.
\end{proof}

\begin{corollary}\label{cor:toggle-invariance}
If $G,H$ are toggle-equivalent, then $q(G;y)=q(H;y)$.
\end{corollary}
\begin{proof}
Since $G,H$ are toggle-equivalent, $M(H)=M(G)*X$ for some feasible $X$ (up to relabeling).
By Lemma~\ref{lem:commute}, $L(M(H))=L(M(G))*\alpha(X)$, so the lifts $G^{\dagger}$ and $H^{\dagger}$ are pivot-equivalent simple graphs.
The vertex-nullity interlace polynomial of a simple graph is invariant under edge pivots and isomorphisms~\cite{ABS2004}, hence under pivot equivalence; therefore $q(G^{\dagger};y)=q(H^{\dagger};y)$.
Theorem~\ref{thm:polynomial} turns this into $y\,q(G;y)=y\,q(H;y)$.
The ring $\mathbb Z[y]$ is an integral domain, so cancellation of $y$ gives $q(G;y)=q(H;y)$.
\end{proof}

\subsection{Quasi-trees and real principal-minor detectors}
The binary principal-matrix theorem for ribbon bouquets states that
\begin{equation}\label{eq:quasitrees}
\QQ(B)=\{X\subseteq E(B):\det_{\F}M(I(B))[X]=1\}
\end{equation}
\cite{Chun2019,DK2026}.
Since the exponent zero contributes the constant polynomial $1$, evaluation at $y=1$ gives
\begin{equation}\label{eq:count}
|\QQ(B)|=q(I(B);1).
\end{equation}

\begin{corollary}\label{cor:quasitree}
For a pseudo-orientable bouquet $B$ and its adjustment $\widehat B$ at any certificate,
$|\QQ(B)|=|\QQ(\widehat B)|$.
\end{corollary}
\begin{proof}
The construction in the converse direction of Theorem~\ref{thm:realization} gives $I(\widehat B)=I(B)^\dagger$.
Therefore
\[
|\QQ(\widehat B)|
=q(I(\widehat B);1)
=q(I(B)^\dagger;1)
=q(I(B);1)
=|\QQ(B)|,
\]
where the first and last equalities are~\eqref{eq:count}, the second substitutes $I(\widehat B)=I(B)^\dagger$, and the third is~\eqref{eq:q} at $y=1$.
\end{proof}

\begin{definition}[Real detectors and PU matrices]\label{def:pu}
A real matrix $A$ is \emph{skew-symmetric} if $A^T=-A$.
It is \emph{principally unimodular}, abbreviated \emph{PU}, if $\det A[X]\in\{0,1,-1\}$ for every principal set $X$.
A \emph{PU skew-symmetric signing} of an alternating binary matrix $H$ is a skew-symmetric integer matrix with entries in $\{0,1,-1\}$ that is PU and reduces entrywise to $H$ modulo $2$.
A \emph{real principal-minor detector} of a symmetric binary matrix $M$ is a real matrix $R$ on the same indices such that
\[
\det_{\mathbb R}R[X]=\det_{\F}M[X]\quad\text{for every }X,
\]
where the binary determinant on the right is interpreted as the integer $0$ or $1$.
\end{definition}

Every odd-order real skew-symmetric determinant vanishes, and every even-order one is nonnegative.
Thus a real skew-symmetric PU matrix has only $0,1$ as principal minors.
The orientable Matrix--Quasi-tree Theorem supplies a PU skew-symmetric signing of the intersection matrix of every orientable bouquet~\cite{MMN2025,DK2026}.

\begin{corollary}[Determinantal counting]\label{cor:determinant}
Let $G$ be stable, $M=M(G)$, and $t=\diag M$.
There exists an integer skew-symmetric matrix $K$ with entries in $\{0,1,-1\}$ such that $R=K+tt^T$ is a real principal-minor detector of $M$.
For every bouquet $B$ with $I(B)\cong G$,
\[
|\QQ(B)|=q(G;1)=\det_{\mathbb R}(I+R).
\]
\end{corollary}
\begin{proof}
Choose an orientable bouquet realizing $G^\dagger$ and a PU skew-symmetric signing $A$ of $L(M)$ supplied by the orientable matrix theorem.
The $p$-column of $L(M)$ is $t$, so $A_{ip}$ is even when $t_i=0$ and is $\pm1$ when $t_i=1$.
Let $D$ be diagonal with $D_{pp}=1$, with $D_{ii}=1$ when $t_i=0$, and with $D_{ii}=A_{ip}$ when $t_i=1$; every diagonal entry of $D$ is $\pm1$.
Then
\[
(DAD)_{ip}=D_{ii}A_{ip}D_{pp}=A_{ip}^2=1\quad(t_i=1),
\]
and $(DAD)_{ip}=0$ when $t_i=0$; similarly the $p$-row becomes $-t^T$, and the $pp$-entry remains $0$.
For every principal set $X$,
\[
\det((DAD)[X])=\det D[X]\cdot\det A[X]\cdot\det D[X]=\det A[X],
\]
because $D[X]$ is diagonal and $\det D[X]^2=1$.
Also $D\equiv I\pmod 2$, so $DAD$ still reduces to $L(M)$ modulo $2$.
Replacing $A$ by $DAD$, we may therefore assume
\[
A=\begin{pmatrix}K&t\\-t^T&0\end{pmatrix},
\]
where $K$ is skew-symmetric, being the top-left block of $A$.

For a real square matrix $C$, its $(i,j)$-cofactor is $(-1)^{i+j}$ times the determinant obtained by deleting row $i$ and column $j$.
The adjugate $\operatorname{adj}(C)$ is the transpose of the cofactor matrix.
Two standard determinant identities are
\[
\det(C+uu^T)=\det C+u^T\operatorname{adj}(C)u
\]
and
\[
\det\begin{pmatrix}C&u\\-u^T&0\end{pmatrix}
=u^T\operatorname{adj}(C)u;
\]
the first is the matrix determinant lemma, and the second follows by expanding the border along its last row.
For an even-order skew-symmetric $C$, its adjugate is skew-symmetric, so the quadratic term is zero.
For an odd-order skew-symmetric $C$, its determinant vanishes, and the second identity applies.

Apply these identities to $C=K[X]$, including the empty case by convention.
If $|X|$ is even, then $\alpha(X)=X$ and $A[\alpha(X)]=K[X]$; the quadratic term vanishes, so
$\det R[X]=\det K[X]=\det A[\alpha(X)]$.
If $|X|$ is odd, then $\alpha(X)=X\cup\{p\}$, the matrix $A[\alpha(X)]$ is the skew border displayed above, and $\det K[X]=0$; hence
$\det R[X]=t[X]^T\operatorname{adj}(K[X])t[X]=\det A[\alpha(X)]$.
In both cases
\[
\det R[X]=\det A[\alpha(X)]\in\{0,1\},
\]
because $A$ is PU.

Also $R\bmod2=M$: indeed $K\bmod2$ is the top-left block of $L(M)$, namely $M+tt^T$, and $R=K+tt^T$.
Hence $\det R[X]\equiv\det M[X]\pmod2$ for every $X$.
Since both prescribed and obtained values belong to $\{0,1\}$, they agree as integers, so $R$ is a real principal-minor detector of $M$.

Finally, expansion by principal minors gives
$\det(I+R)=\sum_X\det R[X]$, and~\eqref{eq:quasitrees} completes the proof.
\end{proof}

The corollary concerns all bouquet realizations of a stable graph, even if a particular realization is not pseudo-orientable.
The reason is that their quasi-tree families are determined by the common binary matrix, whereas a geometric certificate need not be.
The orientable and pseudo-orientable matrix theorems used here are established results~\cite{MMN2025,DK2026}; the corollary is their consequence through the lift.

\begin{example}
For a single looped vertex, $M=[1]$ and $G^\dagger$ is the complete graph $K_2$ on two vertices.
Directly, $q(G;y)=2$ and $q(K_2;y)=2y$.
If $G$ has no loops, $t=0$ and the lift is obtained by adjoining an isolated unlooped vertex, meaning a vertex with no incident edges.
The lift identity then reduces to the usual isolated-vertex factor.
\end{example}

\section{Concluding remarks}
The one-vertex lift simultaneously organizes principal-minor parity, rooted pivot projection, and nullity enumeration.
It yields a geometric existence theorem, a complete finite obstruction construction, and an interlace-polynomial identity.
The exact computation identifies the $22$ minimal toggle classes and supplies checkable certificates for the classification.
Further work may compare these explicit representatives with other descriptions of the same obstruction classes and seek geometric interpretations of individual representatives.

\clearpage
\appendix
\section{The twenty-two representatives}\label{app:representatives}
Each drawing represents one toggle-equivalence class, not every graph in that class.
Vertices are labeled $1,\ldots,n$ and red loops specify the diagonal entries of the adjacency matrix.
In each label, $n$ is the order and $\ell$ is the number of looped vertices.
The drawings are ordinary graph drawings, not chord representations.
Machine-readable matrices and actual double-occurrence words are included in the accompanying data.
\input{enumeration/representatives-figures-en.tex}

\end{document}